# Introduction of Curvilinear Coordinates into Numerical Analysis

Dedicated to Professor Emeritus Yoshiyuki Yamamoto, The University of Tokyo

Hiroshi Isshiki, Institute of Mathematical Analysis, Osaka, Japan,
isshiki@dab.hi-ho.ne.jp

Daisuke Kitazawa, Institute of Industrial Science,
The University of Tokyo, dkita@iis.u-tokyo.ac.jp

**Abstract**

Introduction of curvilinear coordinates might be very convenient in many cases. Theoretically, tensor analysis would be most suited. However, tensor notation can't be used in numerical procedure. For example, the strict discrimination of upper and lower suffices is impossible in the present computer languages. In the present paper, a rather elementary approach more suited to write codes of computer programming is adapted. If we introduce curvilinear coordinates, we would be able to treat in more natural way to handle fixed discontinuity, moving discontinuity and curved boundary. We could generate fine mesh in the neighborhood of fixed and moving discontinuities. Furthermore, if we introduce curvilinear coordinates, we can transform a non-square region into square one and can use a regular mesh. Usually in numerical calculation, a curved boundary is approximated by a jagged or non-smooth boundary. An introduction of curvilinear coordinates could solve many problems in the numerical analyses.

## 1. Introduction

Introduction of curvilinear coordinates might be very convenient in many cases. Theoretically, tensor analysis would be most suited for the purpose. However, tensor notation can't be used conveniently in numerical procedure. For example, the strict discrimination of upper and lower suffices is impossible in the present computer languages[1~3].

In the present paper, a rather classical approach is adapted. The method might be more suited to write codes of computer programming using computer language available now.

In order to treat rapid change of the solution function, we usually use finer mesh to follow the rapid change. If we introduce curvilinear coordinates, we would be able to treat in more natural way to handle fixed discontinuity, moving discontinuity. Namely, If we introduce a curvilinear coordinates, we could generate fine mesh in the neighborhood of fixed and moving discontinuities and increase the accuracy of the solution.

Furthermore, if we introduce curvilinear coordinates, we can transform a non-square region into square one and can use a regular mesh in the transformed coordinates. Furthermore, in numerical calculation, a curved boundary is usually approximated by a jagged or non-smooth boundary.

In sections 3 and 4, fixed and moving discontinuities are discussed, respectively, and, in section 5, transformation of non-square region into square region is discussed. An introduction of curvilinear coordinates could solve many problems in the numerical analyses.

## 2. 1D coordinates transformation

Firstly we discuss 1D case for simplicity. A Curvilinear Coordinates (CC) is introduced through a coordinate transformation:

$$x = x(\xi)\,, \quad \xi = x^{-1}(x)\,. \tag{1}$$

The differentiation between coordinates is summarized below:

$$\frac{d\xi}{dx} = \left(\frac{dx}{d\xi}\right)^{-1}, \tag{2}$$

$$\frac{d^2\xi}{dx^2} = \frac{d}{dx}\left(\frac{dx}{d\xi}\right)^{-1} = \frac{d}{d\xi}\left(\frac{dx}{d\xi}\right)^{-1} \cdot \frac{d\xi}{dx} = -\left(\frac{dx}{d\xi}\right)^{-2} \frac{d^2x}{d\xi^2}\frac{d\xi}{dx} = -\frac{d^2x}{d\xi^2}\left(\frac{dx}{d\xi}\right)^{-3}, \tag{3}$$

$$\frac{d^3\xi}{dx^3} = \frac{d}{dx}\left[-\frac{d^2x}{d\xi^2}\left(\frac{dx}{d\xi}\right)^{-3}\right] = -\frac{d}{d\xi}\left[\frac{d^2x}{d\xi^2}\left(\frac{dx}{d\xi}\right)^{-3}\right] \cdot \frac{d\xi}{dx} = -\frac{d^3x}{d\xi^3}\left(\frac{dx}{d\xi}\right)^{-4} + 3\left(\frac{d^2x}{d\xi^2}\right)^2\left(\frac{dx}{d\xi}\right)^{-5}, \tag{4}$$

$$\begin{aligned}\frac{d^4\xi}{dx^4} &= \frac{d}{dx}\left[-\frac{d^3x}{d\xi^3}\left(\frac{dx}{d\xi}\right)^{-4} + 3\left(\frac{d^2x}{d\xi^2}\right)^2\left(\frac{dx}{d\xi}\right)^{-5}\right] = \frac{d}{d\xi}\left[-\frac{d^3x}{d\xi^3}\left(\frac{dx}{d\xi}\right)^{-4} + 3\left(\frac{d^2x}{d\xi^2}\right)^2\left(\frac{dx}{d\xi}\right)^{-5}\right] \cdot \frac{d\xi}{dx} \\ &= \left[-\frac{d^4x}{d\xi^4}\left(\frac{dx}{d\xi}\right)^{-4} + 4\frac{d^3x}{d\xi^3}\frac{d^2x}{d\xi^2}\left(\frac{dx}{d\xi}\right)^{-5} + 6\frac{d^3x}{d\xi^3}\frac{d^2x}{d\xi^2}\left(\frac{dx}{d\xi}\right)^{-5} - 15\left(\frac{d^2x}{d\xi^2}\right)^2\frac{d^2x}{d\xi^2}\left(\frac{dx}{d\xi}\right)^{-6}\right] \cdot \frac{d\xi}{dx} \\ &= -\frac{d^4x}{d\xi^4}\left(\frac{dx}{d\xi}\right)^{-5} + 10\frac{d^3x}{d\xi^3}\frac{d^2x}{d\xi^2}\left(\frac{dx}{d\xi}\right)^{-6} - 15\left(\frac{d^2x}{d\xi^2}\right)^3\left(\frac{dx}{d\xi}\right)^{-7}.\end{aligned} \tag{5}$$

Coordinates transformation of derivatives of function is summarized as

$$\frac{du}{dx} = \frac{du}{d\xi}\frac{d\xi}{dx}, \tag{6}$$

$$\frac{d^2u}{dx^2} = \frac{d}{dx}\frac{du}{d\xi} \cdot \frac{d\xi}{dx} + \frac{du}{d\xi}\frac{d^2\xi}{dx^2} = \frac{d^2u}{d\xi^2}\left(\frac{d\xi}{dx}\right)^2 + \frac{du}{d\xi}\frac{d^2\xi}{dx^2}, \tag{7}$$

$$\begin{aligned}\frac{d^3u}{dx^3} &= \frac{d}{dx}\frac{d^2u}{d\xi^2} \cdot \left(\frac{d\xi}{dx}\right)^2 + 2\frac{d^2u}{d\xi^2}\frac{d\xi}{dx}\frac{d^2\xi}{dx^2} + \frac{d}{dx}\frac{du}{d\xi} \cdot \frac{d^2\xi}{dx^2} + \frac{du}{d\xi}\frac{d^3\xi}{dx^3} \\ &= \frac{d^3u}{d\xi^3}\left(\frac{d\xi}{dx}\right)^3 + 2\frac{d^2u}{d\xi^2}\frac{d\xi}{dx}\frac{d^2\xi}{dx^2} + \frac{d^2u}{d\xi^2}\frac{d\xi}{dx}\frac{d^2\xi}{dx^2} + \frac{du}{d\xi}\frac{d^3\xi}{dx^3} \\ &= \frac{d^3u}{d\xi^3}\left(\frac{d\xi}{dx}\right)^3 + 3\frac{d^2u}{d\xi^2}\frac{d\xi}{dx}\frac{d^2\xi}{dx^2} + \frac{du}{d\xi}\frac{d^3\xi}{dx^3}.\end{aligned} \tag{8}$$

$$\begin{aligned}\frac{d^4u}{dx^4} &= \frac{d^4u}{d\xi^4}\left(\frac{d\xi}{dx}\right)^4 + 3\frac{d^3u}{d\xi^3}\frac{d^2\xi}{dx^2}\left(\frac{d\xi}{dx}\right)^2 + 3\frac{d^2u}{d\xi^2}\left(\frac{d\xi}{dx}\right)^2\frac{d^2\xi}{dx^2} + 3\frac{d^2u}{d\xi^2}\left(\frac{d^2\xi}{dx^2}\right)^2 + 3\frac{d^2u}{d\xi^2}\frac{d\xi}{dx}\frac{d^3\xi}{dx^3} \\ &\quad + \frac{d^2u}{d\xi^2}\frac{d\xi}{dx}\frac{d^3\xi}{dx^3} + \frac{du}{d\xi}\frac{d^4\xi}{dx^4} \\ &= \frac{d^4u}{d\xi^4}\left(\frac{d\xi}{dx}\right)^4 + 3\frac{d^3u}{d\xi^3}\frac{d^2\xi}{dx^2}\left(\frac{d\xi}{dx}\right)^2 + \frac{d^2u}{d\xi^2}\left[3\left(\frac{d^2\xi}{dx^2}\right)^2 + 4\frac{d\xi}{dx}\frac{d^3\xi}{dx^3}\right] + \frac{du}{d\xi}\frac{d^4\xi}{dx^4}.\end{aligned} \tag{9}$$

For example, we apply the above results to equation of Euler beam:

$$\frac{d^2}{dx^2}\left(EI\frac{d^2u}{dx^2}\right) = f, \tag{10}$$

where $u$, $E$, $I$ and $f$ are the deflection, Young's modulus, second moment of area and external force, respectively. Equation (10) is transformed into

$$\left[\left(\frac{d\xi}{dx}\right)^2\frac{d^2}{d\xi^2} + \frac{d^2\xi}{dx^2}\frac{d}{d\xi}\right]\left(EI\frac{d^2}{dx^2}\left[\left(\frac{d\xi}{dx}\right)^2\frac{d^2u}{d\xi^2} + \frac{d^2\xi}{dx^2}\frac{du}{d\xi}\right]\right) = f. \tag{11}$$

The derivatives $d\xi/dx$ and $d^2\xi/dx^2$ are given by Eqs. (2) and (3), respectively.

## 3. Finite difference solution of 1D burgers equation with a fixed discontinuity

The initial and boundary value problem of Burgers' equation is given by

$$\frac{\partial u}{\partial t}+u\frac{\partial u}{\partial x}=\nu\frac{\partial^2 u}{\partial x^2} \quad \text{in } -L<x<+L \quad \text{for } t>0, \tag{12}$$

$$u=\pm 1 \quad \text{at } x=\mp L \quad \text{for } t>0, \tag{13}$$

$$u=1-\frac{1}{L}(x+L) \quad \text{in } -L<x<+L \quad \text{for } t=0. \tag{14}$$

In this problem, a discontinuity is formed at $x=0$, when time increases.

If we introduce a Curvilinear Coordinates (CC) given by Eq. (1), Eq. (12) is rewritten from Eqs. (6) and (7) as

$$\frac{\partial u}{\partial t}+u\frac{\partial \xi}{\partial x}\frac{\partial u}{\partial \xi}=\nu\left(\left(\frac{\partial \xi}{\partial x}\right)^2\frac{\partial^2 u}{\partial \xi^2}+\frac{\partial^2 \xi}{\partial x^2}\frac{\partial u}{\partial \xi}\right). \tag{15}$$

When $\xi$ is given as an explicit function of $x$, Eq. (15) can be used for numerical calculation.

However, when $x$ is given as an explicit function of $\xi$, we substitute Eqs. (2) and (3) into Eq. (15) and obtain

$$\frac{\partial u}{\partial t}+u\left(\frac{\partial x}{\partial \xi}\right)^{-1}\frac{\partial u}{\partial \xi}=\nu\left(\left(\frac{\partial x}{\partial \xi}\right)^{-2}\frac{\partial^2 u}{\partial \xi^2}-\frac{\partial^2 x}{\partial \xi^2}\left(\frac{\partial x}{\partial \xi}\right)^{-3}\frac{\partial u}{\partial \xi}\right). \tag{16}$$

In the initial and boundary value problem given by Eqs. (12), (13) and (14), the discontinuity is formed at $x=L$ as time increases. Then, we use, for example, the following CC:

$$x=a\left(c\xi+L\xi^{n_p}\right), \quad -1<\xi<+1, \tag{17}$$

where $n_p$ is an odd integer, and $a$ is given by

$$a=\frac{L}{c+L}. \tag{18}$$

Substituting Eq. (17) into Eqs. (2) and (3), we have

$$\frac{\partial x}{\partial \xi}=a(c+3L\xi^2), \quad \frac{\partial^2 x}{\partial \xi^2}=6aL\xi. \tag{19a,b}$$

Substituting Eq. (19) into Eq. (16), we obtain

$$\frac{\partial u}{\partial t}+\frac{1}{a(c+3L\xi^2)}u\frac{\partial u}{\partial \xi}=\nu\left(\frac{1}{a(c+3L\xi^2)}\frac{\partial^2 u}{\partial \xi^2}-\frac{6L\xi}{a^2(c+3L\xi^2)^3}\frac{\partial u}{\partial \xi}\right). \tag{20}$$

From Eq. (3.9), we obtain the following difference equation:

$$\frac{\partial u_i}{\partial t}=-\frac{1}{a(c+3L\xi_i^2)}u_i\frac{u_{i+1}-u_{i-1}}{2d\xi}+\nu\left(\frac{1}{a^2(c+3L\xi_i^2)^2}\frac{u_{i+1}-2u_i+u_{i-1}}{d\xi^2}-\frac{6L\xi_i}{a^2(c+3L\xi_i^2)^3}\frac{u_{i+1}-u_{i-1}}{2d\xi}\right), \tag{21}$$

$$u_i(t+dt)=u_i(t)+\frac{\partial u_i}{\partial t}dt. \tag{22}$$

The boundary and initial conditions are given by

$u\left(x(\xi),t\right)=\pm 1$ at $\xi=\pm 1$ for $t>0$, (23)

$u\left(x(\xi),0\right)=1-\frac{1}{L}\left(x(\xi)+L\right)$ in $-L<\xi<+L$ for $t=0$, (24)

respectively.

In the numerical calculations, parameters are given as $L=8$, $M=80$, $dx=0.2$, $\nu=0.01$ and $dt=0.00125$. Figure 1 shows a coordinates transformation between $x$ and $\xi$ given by Eq. (17). Numerical results are shown in Fig. 2, where CC means Curvilinear Coordinates. If we use the central difference for equal spacing, the solution diverged when t= 10.0. A better result is obtained in CC using central difference as shown in Fig. 2(c).

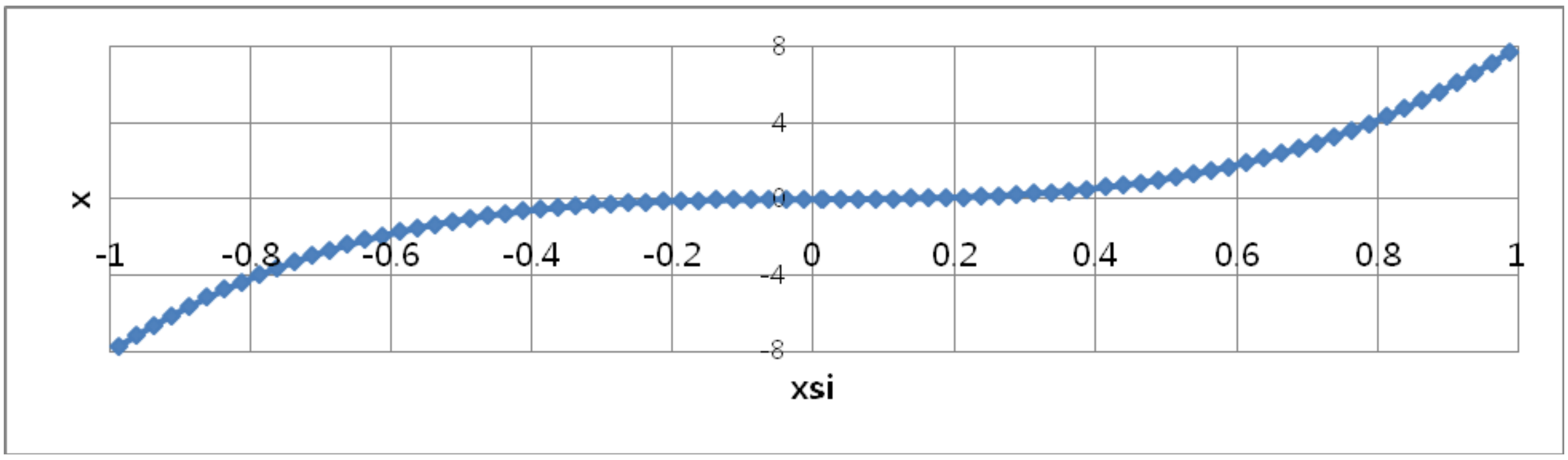

Fig. 1. Coordinates transformation.

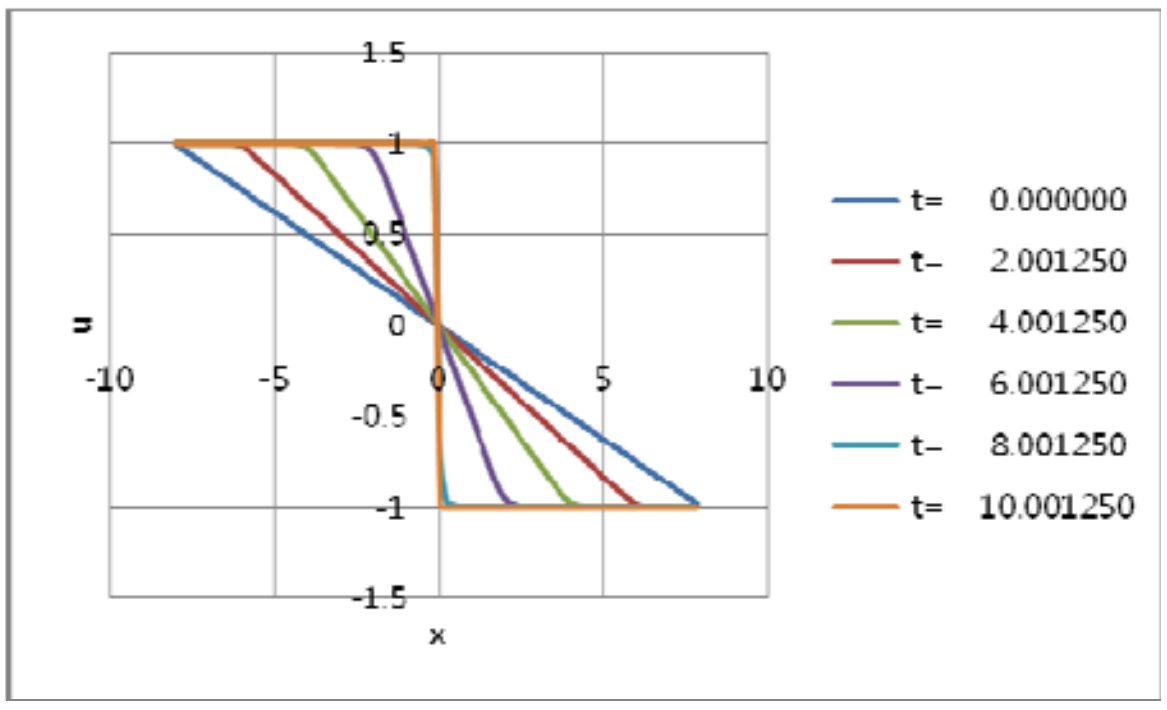

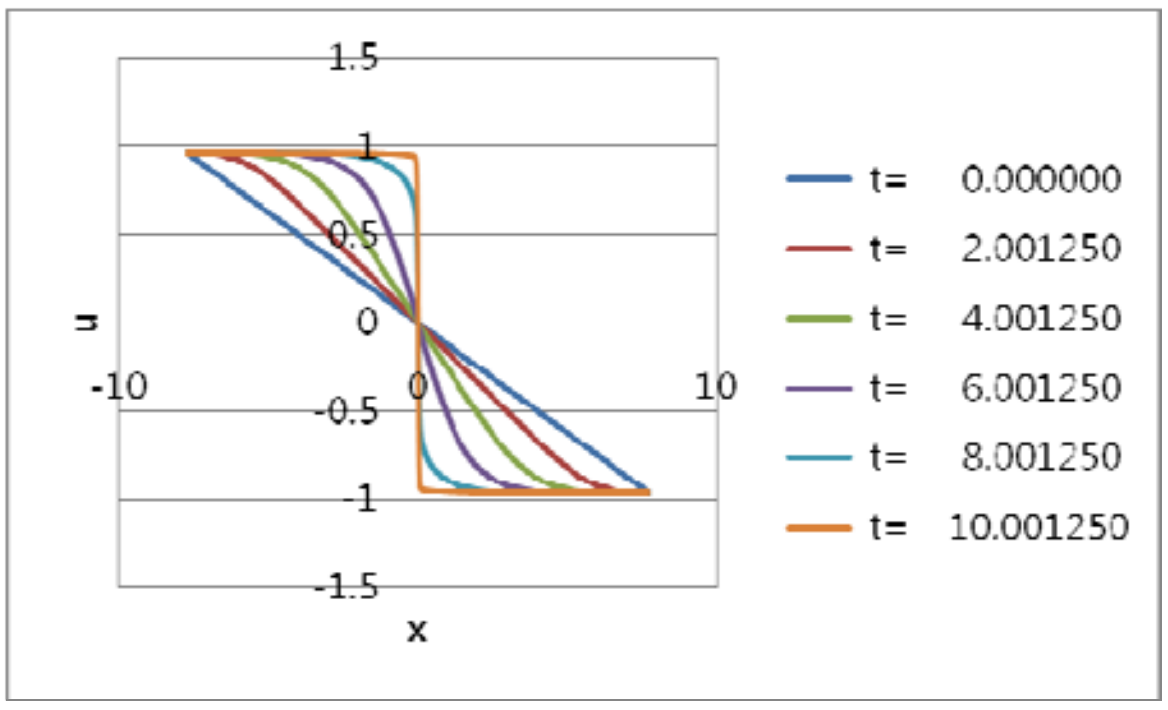

(a) Exact solution (b) Num. sol. in CC (Upwind, $c$=0.2, $n_p$=3)

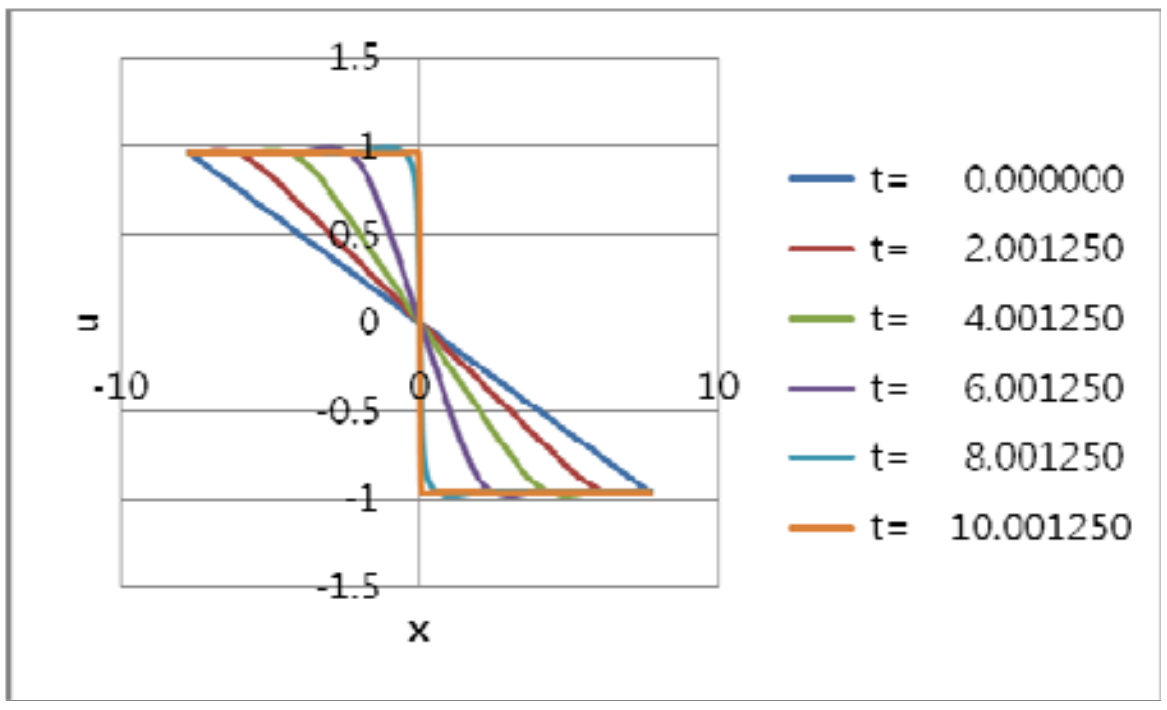

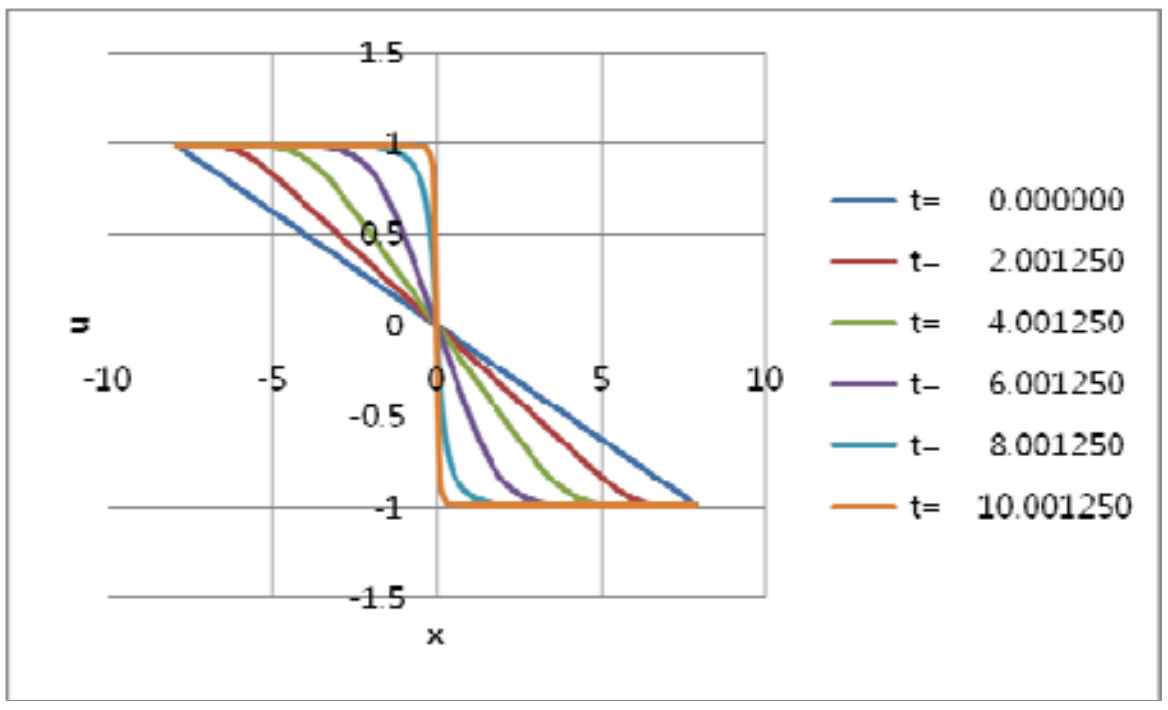

(c) Num. sol. in CC (Central difference, $c$=0.2, $n_p$=3) (d) Num. sol. in Equal spacing (Upwind, $c$=0, $n_p$=1)

Fig. 2. Finite Difference solution of 1D burgers equation with a fixed discontinuity ($M$=80)

**4. Finite Difference solution of 1D convective diffusion equation with a moving discontinuity**

We consider the following convective diffusion problem:

$$\frac{\partial u}{\partial t} + U\frac{\partial u}{\partial x} = \nu \frac{\partial^2 u}{\partial x^2} \quad \text{in} \quad -L < x < +L \quad \text{for} \quad t > 0, \tag{25}$$

$$u = 1 \quad \text{at} \quad x = -L \quad \text{and} \quad u = 0 \quad \text{at} \quad x = -L \quad \text{for} \quad t > 0, \tag{26}$$

$$u = 1 \quad \text{in} \quad -L < x \le -L/4 \quad \text{and} \quad u = 0 \quad \text{in} \quad -L/4 < x < L \quad \text{for} \quad t = 0. \tag{27}$$

In this initial and boundary value problem, the discontinuity moves with velocity U.

We introduce Curvilinear Coordinates (CC) $(\xi, \tau)$ as

$$x = x(\xi, \tau), \quad t = \tau, \tag{28a, b}$$

$$\xi = \xi(x,t), \quad \tau = t. \tag{29a, b}$$

From Eqs. (28) and (29), we have

$$\frac{\partial}{\partial x} = \frac{\partial \xi}{\partial x}\frac{\partial}{\partial \xi} + \frac{\partial \tau}{\partial x}\frac{\partial}{\partial \tau} = \frac{\partial \xi}{\partial x}\frac{\partial}{\partial \xi}, \quad \frac{\partial}{\partial t} = \frac{\partial \xi}{\partial t}\frac{\partial}{\partial \xi} + \frac{\partial \tau}{\partial t}\frac{\partial}{\partial \tau} = \frac{\partial \xi}{\partial t}\frac{\partial}{\partial \xi} + \frac{\partial}{\partial \tau}, \tag{30a, b}$$

$$\frac{\partial}{\partial \xi} = \frac{\partial x}{\partial \xi}\frac{\partial}{\partial x} + \frac{\partial t}{\partial \xi}\frac{\partial}{\partial t} = \frac{\partial x}{\partial \xi}\frac{\partial}{\partial x}, \quad \frac{\partial}{\partial \tau} = \frac{\partial x}{\partial \tau}\frac{\partial}{\partial x} + \frac{\partial t}{\partial \tau}\frac{\partial}{\partial t} = \frac{\partial x}{\partial \tau}\frac{\partial}{\partial x} + \frac{\partial}{\partial t}. \tag{31a, b}$$

From Eq. (30a), we obtain

$$1 = \frac{\partial \xi}{\partial x}\frac{\partial x}{\partial \xi} \quad \text{or} \quad \frac{\partial \xi}{\partial x} = \left(\frac{\partial x}{\partial \xi}\right)^{-1}. \tag{31}$$

From Eq. (31b), we obtain

$$0 = \frac{\partial x}{\partial \tau}\frac{\partial \xi}{\partial x} + \frac{\partial \xi}{\partial t} \quad \text{or} \quad \frac{\partial \xi}{\partial t} = -\frac{\partial x}{\partial \tau}\frac{\partial \xi}{\partial x} = -\frac{\partial x}{\partial \tau}\left(\frac{\partial x}{\partial \xi}\right)^{-1}. \tag{32}$$

From Eq. (30), the derivatives of $u$ can be written as

$$\frac{\partial u}{\partial x} = \frac{\partial \xi}{\partial x}\frac{\partial u}{\partial \xi} + \frac{\partial \tau}{\partial x}\frac{\partial u}{\partial \tau} = \frac{\partial \xi}{\partial x}\frac{\partial u}{\partial \xi}, \quad \frac{\partial u}{\partial t} = \frac{\partial \xi}{\partial t}\frac{\partial u}{\partial \xi} + \frac{\partial \tau}{\partial t}\frac{\partial u}{\partial \tau} = \frac{\partial \xi}{\partial t}\frac{\partial u}{\partial \xi} + \frac{\partial u}{\partial \tau}, \tag{33a, b}$$

$$\begin{aligned}\frac{\partial^2 u}{\partial x^2} &= \frac{\partial}{\partial x}\left(\frac{\partial \xi}{\partial x}\frac{\partial u}{\partial \xi}\right) = \frac{\partial \xi}{\partial x}\frac{\partial}{\partial \xi}\left(\frac{\partial \xi}{\partial x}\frac{\partial u}{\partial \xi}\right) + \frac{\partial \tau}{\partial x}\frac{\partial}{\partial \tau}\left(\frac{\partial \xi}{\partial x}\frac{\partial u}{\partial \xi}\right) = \frac{\partial \xi}{\partial x}\frac{\partial}{\partial \xi}\left(\frac{\partial \xi}{\partial x}\right)\cdot\frac{\partial u}{\partial \xi} + \frac{\partial \xi}{\partial x}\frac{\partial \xi}{\partial x}\frac{\partial}{\partial \xi}\left(\frac{\partial u}{\partial \xi}\right) \\ &= \frac{\partial \xi}{\partial x}\frac{\partial}{\partial \xi}\left(\frac{\partial \xi}{\partial x}\right)\cdot\frac{\partial u}{\partial \xi} + \left(\frac{\partial \xi}{\partial x}\right)^2\frac{\partial^2 u}{\partial \xi^2}.\end{aligned} \tag{34}$$

Substituting Eqs. (31) and (32) into Eqs. (33) and (34), we have

$$\frac{\partial u}{\partial x} = \left(\frac{\partial x}{\partial \xi}\right)^{-1}\frac{\partial u}{\partial \xi}, \quad \frac{\partial u}{\partial t} = -\frac{\partial x}{\partial \tau}\left(\frac{\partial x}{\partial \xi}\right)^{-1}\frac{\partial u}{\partial \xi} + \frac{\partial u}{\partial \tau}, \tag{35a, b}$$

$$\frac{\partial^2 u}{\partial x^2} = \left(\frac{\partial x}{\partial \xi}\right)^{-1}\frac{\partial}{\partial \xi}\left(\frac{\partial x}{\partial \xi}\right)^{-1}\frac{\partial u}{\partial \xi} + \left(\frac{\partial x}{\partial \xi}\right)^{-2}\frac{\partial^2 u}{\partial \xi^2} = -\left(\frac{\partial x}{\partial \xi}\right)^{-3}\frac{\partial^2 x}{\partial \xi^2}\frac{\partial u}{\partial \xi} + \left(\frac{\partial x}{\partial \xi}\right)^{-2}\frac{\partial^2 u}{\partial \xi^2}. \tag{36}$$

Substituting Eqs. (35) and (36) into Eq. (25), we obbtain

$$-\frac{\partial x}{\partial \tau}\left(\frac{\partial x}{\partial \xi}\right)^{-1}\frac{\partial u}{\partial \xi} + \frac{\partial u}{\partial \tau} + U\left(\frac{\partial x}{\partial \xi}\right)^{-1}\frac{\partial u}{\partial \xi} = \nu\left[-\left(\frac{\partial x}{\partial \xi}\right)^{-3}\frac{\partial^2 x}{\partial \xi^2}\frac{\partial u}{\partial \xi} + \left(\frac{\partial x}{\partial \xi}\right)^{-2}\frac{\partial^2 u}{\partial \xi^2}\right]. \tag{37}$$

Rewriting Eq. (37), the convective diffusion equation (25) in $(\xi, \tau)$ system can be written as

$$\frac{\partial u}{\partial \tau} + \left(U - \frac{\partial x}{\partial \tau}\right)\left(\frac{\partial x}{\partial \xi}\right)^{-1}\frac{\partial u}{\partial \xi} = \nu\left(\frac{\partial x}{\partial \xi}\right)^{-2}\left[\frac{\partial^2 u}{\partial \xi^2} - \left(\frac{\partial x}{\partial \xi}\right)^{-1}\frac{\partial^2 x}{\partial \xi^2}\frac{\partial u}{\partial \xi}\right]. \tag{38}$$

The image of the coordinates transformation is illustrated in Fig.3.

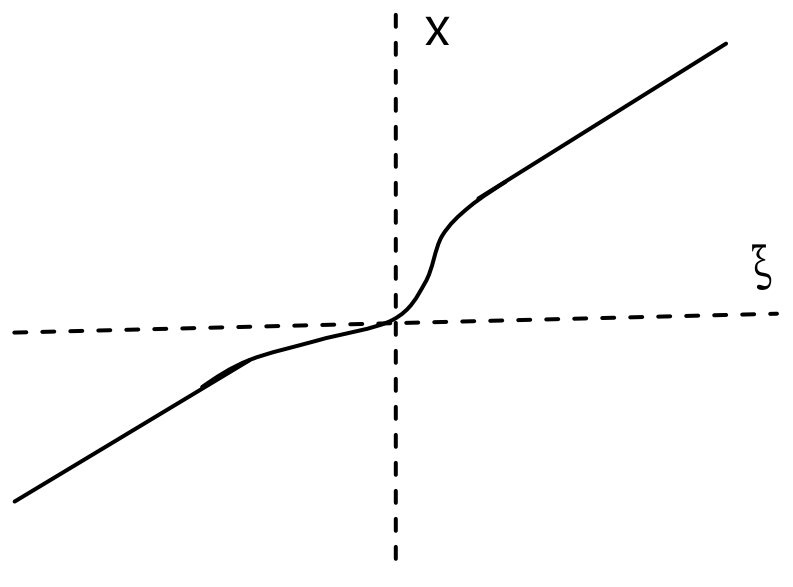

(a) $x=x(\xi,\tau)$

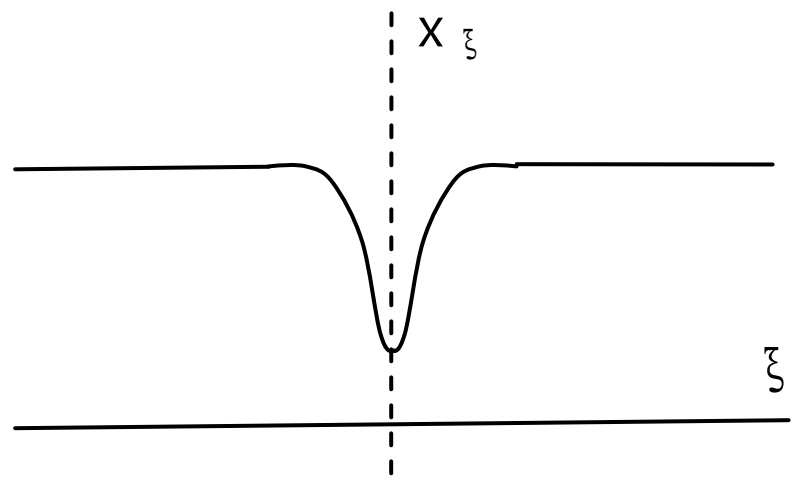

(b) $x_\xi=x_\xi(\xi,\tau)$

Fig. 3. Curvilinear Coordinates (CC)

Let $\xi_0(\tau)$ be the position of the discontinuity at time $\tau$. When $b$ in the following expression is big, the coordinates transformation could be given by

$$\frac{\partial x}{\partial \xi} = L\left(1 + \frac{\sqrt{\pi}\,h}{2b} - he^{-(b(\xi-\xi_0(\tau)))^2}\right), \tag{39}$$

since

$$\int_{-1}^{+1} \frac{\partial x}{\partial \xi} d\xi = 2L \quad \text{and} \quad \int_{-1}^{+1} he^{-(b(\xi-\xi_0(\tau)))^2} d\xi \approx \int_{-\infty}^{+\infty} he^{-(b(\xi-\xi_0(\tau)))^2} d\xi = \frac{h}{b}\int_{-\infty}^{+\infty} e^{-p^2} dp = \frac{\sqrt{\pi}\,h}{b}. \tag{40}$$

Integrating Eq. (39), we obtain

$$\begin{aligned} x &= L\left(\left(1 + \frac{\sqrt{\pi}\,h}{2b}\right)\xi - \int_0^{\xi} he^{-(b(\xi-\xi_0(\tau)))^2} d\xi + C(\tau)\right) = L\left(\left(1 + \frac{\sqrt{\pi}\,h}{2b}\right)\xi - \frac{h}{b}\int_{-b\xi_0(\tau)}^{b(\xi-\xi_0(\tau))} e^{-p^2} dp + C(\tau)\right) \\ &= L\left(\left(1 + \frac{\sqrt{\pi}\,h}{2b}\right)\xi - \frac{h}{b}\int_{-b\xi_0(\tau)}^{0} e^{-p^2} dp - \frac{h}{b}\int_{0}^{b(\xi-\xi_0(\tau))} e^{-p^2} dp + C(\tau)\right) \\ &= L\left(\left(1 + \frac{\sqrt{\pi}\,h}{2b}\right)\xi - \frac{\sqrt{\pi}\,h}{2b}\operatorname{erf}\left(b\xi_0(\tau)\right) - \frac{\sqrt{\pi}\,h}{2b}\operatorname{erf}\left(b(\xi-\xi_0(\tau))\right) + C(\tau)\right), \end{aligned} \tag{41}$$

where $\operatorname{erf}(x)$ is the error function, and $C(\tau)$ is an arbitrary function of $\tau$. We could choose $C(\tau)$ as

$$C(\tau) = \frac{\sqrt{\pi}\,h}{2b}\operatorname{erf}\left(b\xi_0(\tau)\right). \tag{42}$$

Hence, we have

$$x = L\left(\left(1 + \frac{\sqrt{\pi}\,h}{2b}\right)\xi - \frac{\sqrt{\pi}\,h}{2b}\operatorname{erf}\left(b(\xi-\xi_0(\tau))\right)\right). \tag{43}$$

When $x$ is given, $\xi$ is obtained by solving Eq. (43) by using Newton-Raphson method. The correction of $\xi$ is given by

$$d\xi = \left(\frac{x}{L} - \left(1 + \frac{\sqrt{\pi}\,h}{2b}\right)\xi + \frac{\sqrt{\pi}\,h}{2b}\operatorname{erf}\left(b(\xi-\xi_0(\tau))\right)\right)\left(\left(1 + \frac{\sqrt{\pi}\,h}{2b}\right) - h\exp\left(-(b(\xi-\xi_0(\tau)))^2\right)\right)^{-1}. \tag{44}$$

From Eq. (43), we have

$$\frac{\partial x}{\partial \tau} = L\frac{h}{b}\exp\left(-(b(\xi-\xi_0(\tau)))^2\right)b\frac{d\xi_0(\tau)}{d\tau}. \tag{45}$$

From Eq. (43), we obtain

$$\frac{\partial^2 x}{\partial \xi^2} = \frac{\partial}{\partial \xi}L\left(1 + \frac{\sqrt{\pi}\,h}{2b} - he^{-(b(\xi-\xi_0(\tau)))^2}\right) = -Lh\frac{\partial e^{-(b(\xi-\xi_0(\tau)))^2}}{\partial \xi} = 2Lb^2he^{-(b(\xi-\xi_0(\tau)))^2}(\xi-\xi_0(\tau)). \tag{46}$$

The position of the discontinuity $\xi_0(\tau)$ is obtained as follows. Let $x_0 = -L/2$ be the initial position of the discontinuity. Then, $\xi_0(\tau)$ satisfies

$$x_0 + Ut = L\left(\left(1 + \frac{\sqrt{\pi}\,h}{2b}\right)\xi_0(\tau) - \frac{\sqrt{\pi}\,h}{2b}\operatorname{erf}\left(b(\xi_0(\tau)-\xi_0(\tau))\right)\right) = L\left(1 + \frac{\sqrt{\pi}\,h}{2b}\right)\xi_0(\tau). \tag{47}$$

Since $t = \tau$, we rewrite Eq. (47):

$$\xi_0(\tau) = \frac{1}{L}\left(1 + \frac{\sqrt{\pi}\,h}{2b}\right)^{-1}(x_0 + U\tau). \tag{48}$$

From Eq. (48), we have

$$\frac{d\xi_0(\tau)}{d\tau} = \frac{U}{L}\left(1 + \frac{\sqrt{\pi}\, h}{2b}\right)^{-1} . \tag{49}$$

In the numerical calculations, parameters are given as $L = 8$, $M = 80$ or $160$, $dx = 0.2$, $\nu = 0.01$ and $dt = 0.0025$. Numerical results are shown in Figs. 4 and 5, where CC means Curvilinear Coordinates. A much better result is obtained in CC.

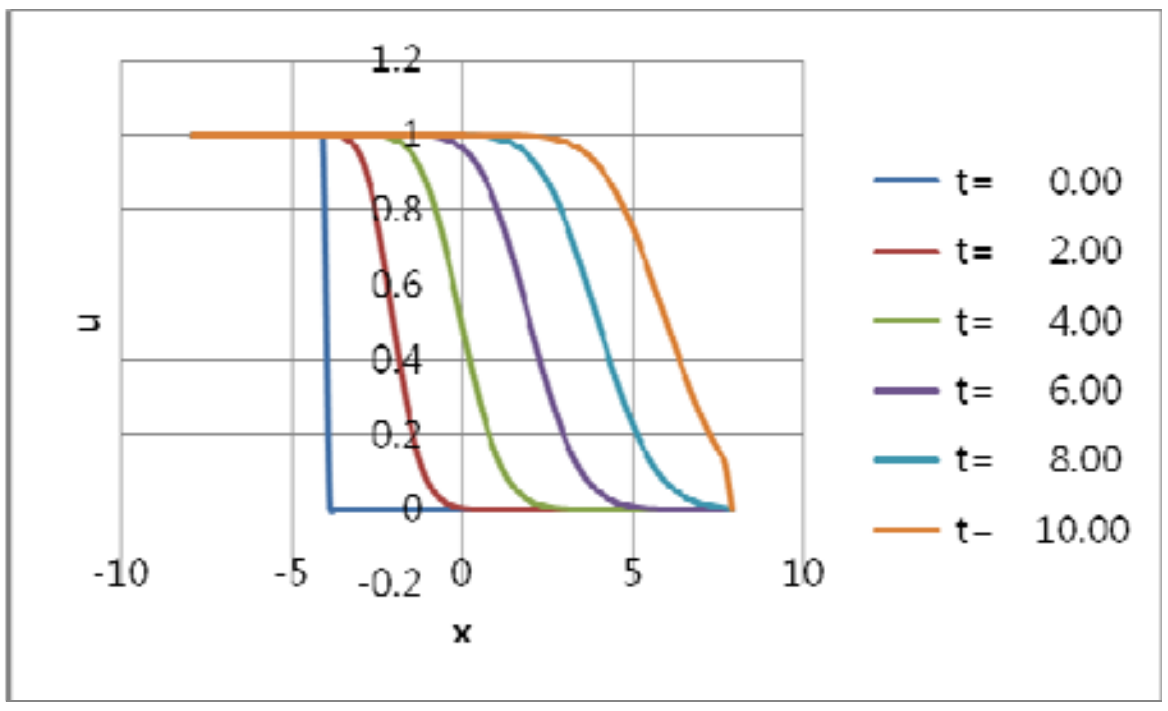

(a) Equal spacing in $x$ (Upwind, $M$=80, $h$=0)

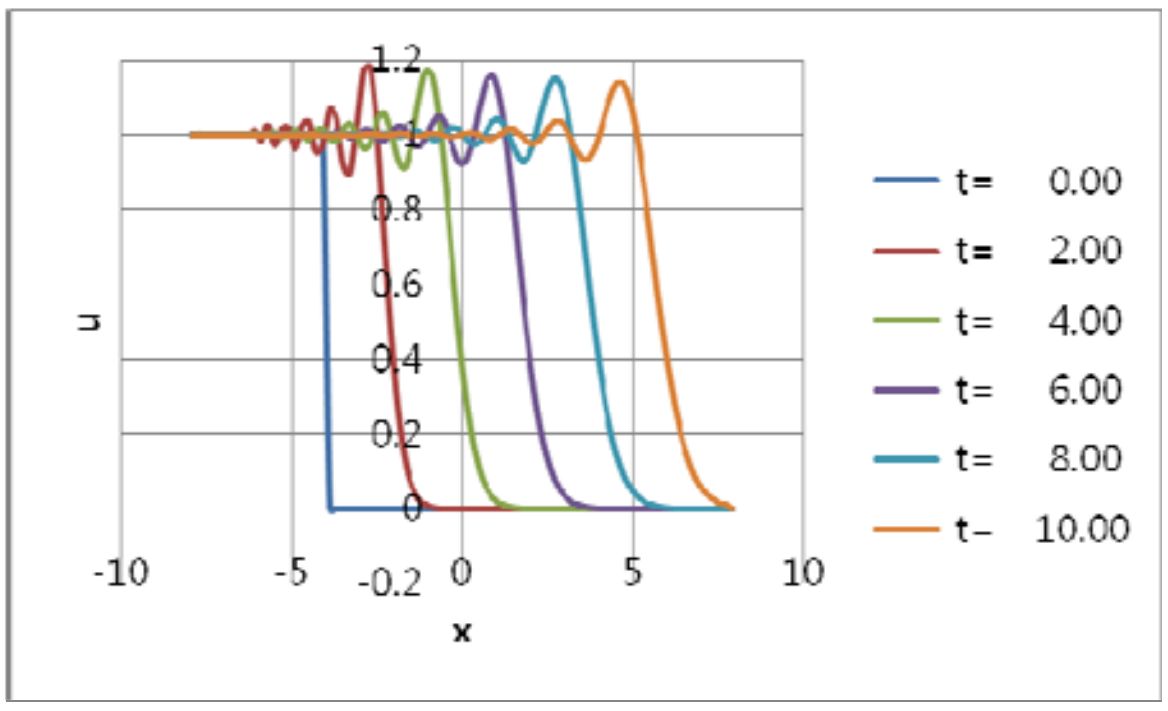

(b) Equal spacing in $x$ (Central difference, $M$=80, $h$=0)

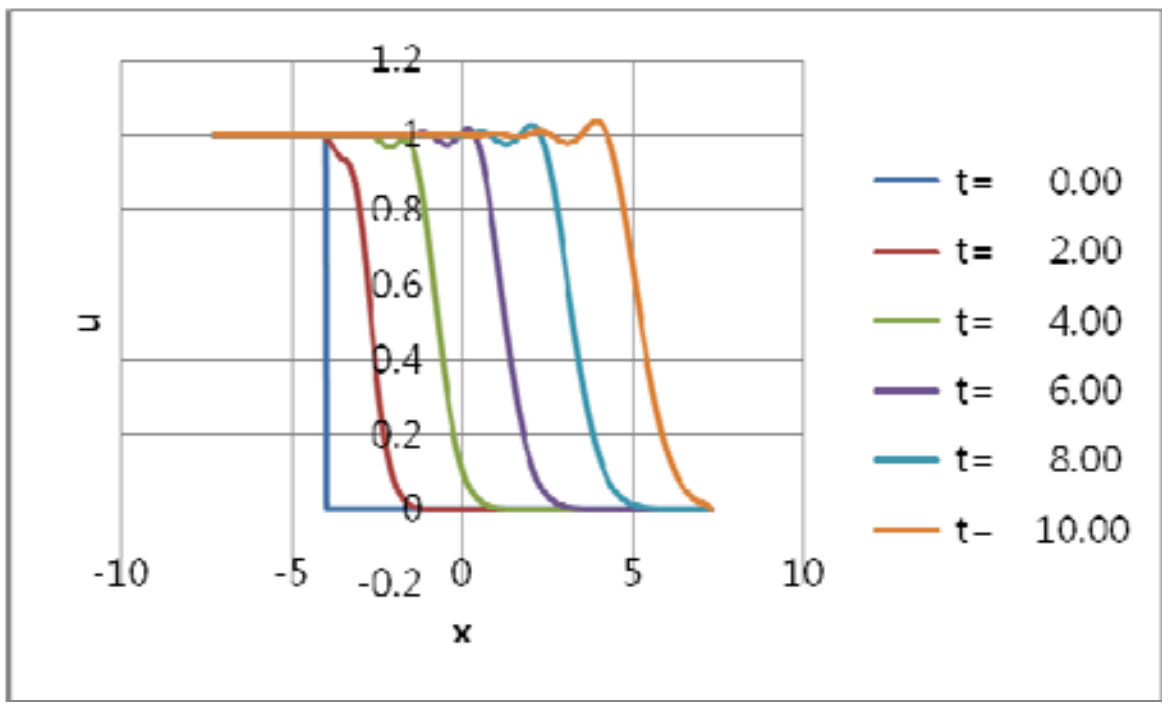

(c) Equal spacing in $\xi$ or CC (Central difference, $M$=80, $h$=0.9, $b$=10)

Fig. 4. Transition of discontinuity ($M$=80)

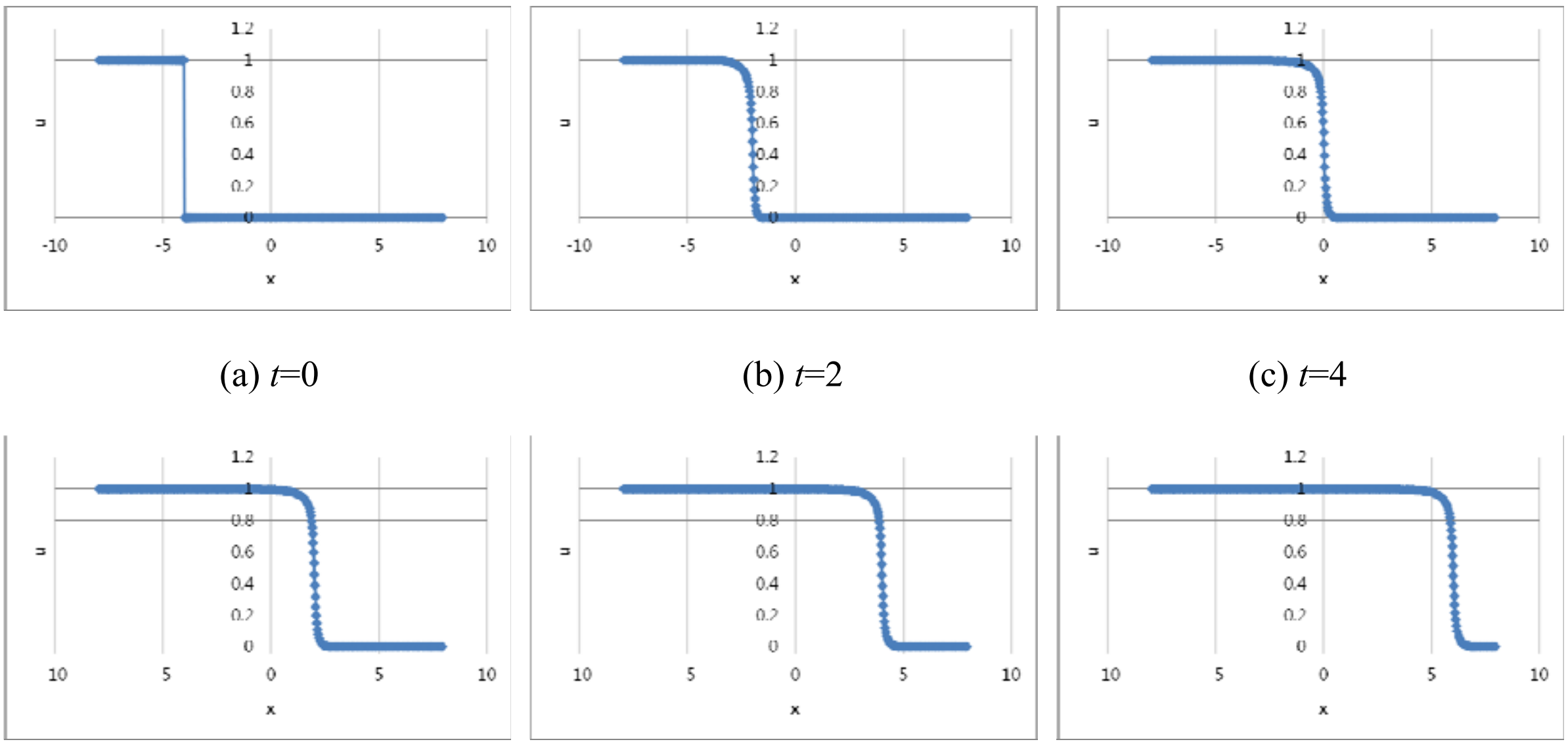

(a) $t$=0 (b) $t$=2 (c) $t$=4

(d) $t$=6 (e) $t$=8 (f) $t$=10

Fig. 5. Transition of discontinuity (Equal spacing in $\xi$ or CC (Central difference, $M$=160, $h$=0.99, $b$=5)

## 5. Potential flow around a circle

We now apply Curvilinear Coordinates (CC) to solve potential flow around a circular cylinder in uniform flow. In the conventional numerical procedure, the circular boundary is approximated by a jagged or non-smooth boundary. However, if we use CC, we could handle the curved boundary more reasonably.

The steady boundary value problem of a potential flow around a circular cylinder in a uniform flow with velocity $U$ is given by

$$\Phi = Ux + \phi\,, \tag{50}$$

$$\frac{\partial^2\phi}{\partial x^2} + \frac{\partial^2\phi}{\partial y^2} = 0 \quad \text{in} \quad a < \sqrt{x^2+y^2} < \infty\,, \quad -\pi < \tan^{-1}\frac{y}{x} < +\pi\,, \tag{51}$$

$$\frac{\partial\phi}{\partial x}n_x + \frac{\partial\phi}{\partial y}n_y = -Un_x \quad \text{on} \quad \sqrt{x^2+y^2} = a\,, \quad -\pi \le \tan^{-1}\frac{y}{x} \le +\pi\,, \tag{52a}$$

$$\phi \to 0 \quad \text{as} \quad \sqrt{x^2+y^2} \to +\infty\,. \tag{52a}$$

where $(x, y)$ is Cartesian coordinates with the origin at the center of the circle, $\Phi$ and $\phi$ are velocity potentials, and $(n_x, n_y)$ is unit normal vector on the circle inward to fluid.

We define a CC $(\xi, \eta)$ as

$$x = x(\xi,\eta)\,, \qquad y = y(\xi,\eta)\,, \tag{53}$$

$$\xi = \xi(x,y)\,, \qquad \eta = \eta(x,y)\,. \tag{54}$$

If we apply

$$\begin{bmatrix} \dfrac{\partial}{\partial x} \\ \dfrac{\partial}{\partial y} \end{bmatrix} = \begin{bmatrix} \dfrac{\partial \xi}{\partial x} & \dfrac{\partial \eta}{\partial x} \\ \dfrac{\partial \xi}{\partial y} & \dfrac{\partial \eta}{\partial y} \end{bmatrix} \begin{bmatrix} \dfrac{\partial}{\partial \xi} \\ \dfrac{\partial}{\partial \eta} \end{bmatrix} \tag{55}$$

and

$$\begin{bmatrix} \dfrac{\partial}{\partial \xi} \\ \dfrac{\partial}{\partial \eta} \end{bmatrix} = \begin{bmatrix} \dfrac{\partial x}{\partial \xi} & \dfrac{\partial y}{\partial \xi} \\ \dfrac{\partial x}{\partial \eta} & \dfrac{\partial y}{\partial \eta} \end{bmatrix} \begin{bmatrix} \dfrac{\partial}{\partial x} \\ \dfrac{\partial}{\partial y} \end{bmatrix}, \tag{56}$$

we have

$$\begin{bmatrix} 1 & 0 \\ 0 & 1 \end{bmatrix} = \begin{bmatrix} \dfrac{\partial x}{\partial \xi} & \dfrac{\partial x}{\partial \eta} \\ \dfrac{\partial y}{\partial \xi} & \dfrac{\partial y}{\partial \eta} \end{bmatrix} \begin{bmatrix} \dfrac{\partial \xi}{\partial x} & \dfrac{\partial \xi}{\partial y} \\ \dfrac{\partial \eta}{\partial x} & \dfrac{\partial \eta}{\partial y} \end{bmatrix}. \tag{57}$$

Equation (57) is obtained by taking the transpose of a matrix expression derived by substituting $x$ and $y$ into Eq. (55). From Eq. (57), we obtain

$$\begin{bmatrix} \dfrac{\partial \xi}{\partial x} & \dfrac{\partial \xi}{\partial y} \\ \dfrac{\partial \eta}{\partial x} & \dfrac{\partial \eta}{\partial y} \end{bmatrix} = \begin{bmatrix} \dfrac{\partial x}{\partial \xi} & \dfrac{\partial x}{\partial \eta} \\ \dfrac{\partial y}{\partial \xi} & \dfrac{\partial y}{\partial \eta} \end{bmatrix}^{-1} \tag{58}$$

or

$$\begin{aligned} \frac{\partial \xi}{\partial x} &= \left(g^{-1}\right)_{11} = \left(\frac{\partial x}{\partial \xi}\frac{\partial y}{\partial \eta} - \frac{\partial x}{\partial \eta}\frac{\partial y}{\partial \xi}\right)^{-1}\frac{\partial y}{\partial \eta} = \frac{1}{|g|}\frac{\partial y}{\partial \eta} \\ \frac{\partial \eta}{\partial x} &= \left(g^{-1}\right)_{21} = -\left(\frac{\partial x}{\partial \xi}\frac{\partial y}{\partial \eta} - \frac{\partial x}{\partial \eta}\frac{\partial y}{\partial \xi}\right)^{-1}\frac{\partial y}{\partial \xi} = -\frac{1}{|g|}\frac{\partial y}{\partial \xi} \end{aligned}, \tag{59a}$$

$$\begin{aligned} \frac{\partial \xi}{\partial y} &= \left(g^{-1}\right)_{12} = -\left(\frac{\partial x}{\partial \xi}\frac{\partial y}{\partial \eta} - \frac{\partial x}{\partial \eta}\frac{\partial y}{\partial \xi}\right)^{-1}\frac{\partial x}{\partial \eta} = -\frac{1}{|g|}\frac{\partial x}{\partial \eta} \\ \frac{\partial \eta}{\partial y} &= \left(g^{-1}\right)_{22} = \left(\frac{\partial x}{\partial \xi}\frac{\partial y}{\partial \eta} - \frac{\partial x}{\partial \eta}\frac{\partial y}{\partial \xi}\right)^{-1}\frac{\partial x}{\partial \xi} = \frac{1}{|g|}\frac{\partial x}{\partial \xi} \end{aligned}, \tag{59b}$$

where

$$g = \begin{bmatrix} \dfrac{\partial x}{\partial \xi} & \dfrac{\partial x}{\partial \eta} \\ \dfrac{\partial y}{\partial \xi} & \dfrac{\partial y}{\partial \eta} \end{bmatrix}, \quad |g| = \frac{\partial x}{\partial \xi}\frac{\partial y}{\partial \eta} - \frac{\partial x}{\partial \eta}\frac{\partial y}{\partial \xi}. \tag{60a,b}$$

Similarly, from Eq. (57), we obtain, we obtain

$$\begin{bmatrix} \dfrac{\partial x}{\partial \xi} & \dfrac{\partial x}{\partial \eta} \\ \dfrac{\partial y}{\partial \xi} & \dfrac{\partial y}{\partial \eta} \end{bmatrix} = \begin{bmatrix} \dfrac{\partial \xi}{\partial x} & \dfrac{\partial \xi}{\partial y} \\ \dfrac{\partial \eta}{\partial x} & \dfrac{\partial \eta}{\partial y} \end{bmatrix}^{-1} \tag{61}$$

or

$$\begin{aligned} \frac{\partial x}{\partial \xi} &= G_{11} = \left(\frac{\partial \xi}{\partial x}\frac{\partial \eta}{\partial y} - \frac{\partial \xi}{\partial y}\frac{\partial \eta}{\partial x}\right)^{-1}\frac{\partial \eta}{\partial y} = \frac{1}{|G|}\frac{\partial \eta}{\partial y} \\ \frac{\partial y}{\partial \xi} &= G_{21} = -\left(\frac{\partial \xi}{\partial x}\frac{\partial \eta}{\partial y} - \frac{\partial \xi}{\partial y}\frac{\partial \eta}{\partial x}\right)^{-1}\frac{\partial \eta}{\partial x} = -\frac{1}{|G|}\frac{\partial \eta}{\partial x} \end{aligned}, \tag{62a}$$

$$\begin{aligned} \frac{\partial x}{\partial \eta} &= G_{12} = -\left(\frac{\partial \xi}{\partial x}\frac{\partial \eta}{\partial y} - \frac{\partial \xi}{\partial y}\frac{\partial \eta}{\partial x}\right)^{-1}\frac{\partial \xi}{\partial y} = -\frac{1}{|G|}\frac{\partial \xi}{\partial y} \\ \frac{\partial y}{\partial \eta} &= G_{22} = \left(\frac{\partial \xi}{\partial x}\frac{\partial \eta}{\partial y} - \frac{\partial \xi}{\partial y}\frac{\partial \eta}{\partial x}\right)^{-1}\frac{\partial \xi}{\partial x} = \frac{1}{|G|}\frac{\partial \xi}{\partial x} \end{aligned}, \tag{62b}$$

where

$$G = \begin{bmatrix} \dfrac{\partial \xi}{\partial x} & \dfrac{\partial \xi}{\partial y} \\ \dfrac{\partial \eta}{\partial x} & \dfrac{\partial \eta}{\partial y} \end{bmatrix}, \quad |G| = \frac{\partial \xi}{\partial x}\frac{\partial \eta}{\partial y} - \frac{\partial \eta}{\partial x}\frac{\partial \xi}{\partial y}. \tag{63a,b}$$

With respect to matrices $g$ and $G$, we can derive the following properties:

$$\frac{\partial |g|}{\partial \xi} = \frac{\partial}{\partial \xi}\left( \frac{\partial x}{\partial \xi}\frac{\partial y}{\partial \eta} - \frac{\partial x}{\partial \eta}\frac{\partial y}{\partial \xi} \right) = \frac{\partial^2 x}{\partial \xi^2}\frac{\partial y}{\partial \eta} + \frac{\partial x}{\partial \xi}\frac{\partial^2 y}{\partial \xi \partial \eta} - \frac{\partial^2 x}{\partial \xi \partial \eta}\frac{\partial y}{\partial \xi} - \frac{\partial x}{\partial \eta}\frac{\partial^2 y}{\partial \xi^2}, \tag{64a}$$

$$\frac{\partial |g|}{\partial \eta} = \frac{\partial}{\partial \eta}\left( \frac{\partial x}{\partial \xi}\frac{\partial y}{\partial \eta} - \frac{\partial x}{\partial \eta}\frac{\partial y}{\partial \xi} \right) = \frac{\partial^2 x}{\partial \xi \partial \eta}\frac{\partial y}{\partial \eta} + \frac{\partial x}{\partial \xi}\frac{\partial^2 y}{\partial \eta^2} - \frac{\partial^2 x}{\partial \eta^2}\frac{\partial y}{\partial \xi} - \frac{\partial x}{\partial \eta}\frac{\partial^2 y}{\partial \xi \partial \eta} \tag{64b}$$

and

$$\frac{\partial |G|}{\partial x} = \frac{\partial}{\partial x}\left( \frac{\partial \xi}{\partial x}\frac{\partial \eta}{\partial y} - \frac{\partial \xi}{\partial y}\frac{\partial \eta}{\partial x} \right) = \frac{\partial^2 \xi}{\partial x^2}\frac{\partial \eta}{\partial y} + \frac{\partial \xi}{\partial x}\frac{\partial^2 \eta}{\partial x \partial y} - \frac{\partial^2 \xi}{\partial x \partial y}\frac{\partial \eta}{\partial x} - \frac{\partial \xi}{\partial y}\frac{\partial^2 \eta}{\partial x^2}, \tag{65a}$$

$$\frac{\partial |G|}{\partial y} = \frac{\partial}{\partial y}\left( \frac{\partial \xi}{\partial x}\frac{\partial \eta}{\partial y} - \frac{\partial \xi}{\partial y}\frac{\partial \eta}{\partial x} \right) = \frac{\partial^2 \xi}{\partial x \partial y}\frac{\partial \eta}{\partial y} + \frac{\partial \xi}{\partial x}\frac{\partial^2 \eta}{\partial y^2} - \frac{\partial^2 \xi}{\partial y^2}\frac{\partial \eta}{\partial x} - \frac{\partial \xi}{\partial y}\frac{\partial^2 \eta}{\partial x \partial y}. \tag{65b}$$

From Eq. (57), we have

$$gG = I. \tag{66}$$

When $(x, y)$ is expressed explicitly as a function of $(\xi, \eta)$, matrix $g$ and Eq. (59) become important. On the contrary, if $(\xi, \eta)$ is expressed explicitly as a function of $(x, y)$, matrix $G$ and Eq. (62) become important.

When $(x, y)$ is expressed explicitly as a function of $(\xi, \eta)$, the second derivatives of $(\xi, \eta)$ with respect to $(x, y)$ are given by

$$\begin{aligned}
\frac{\partial^2 \xi}{\partial x^2} &= \frac{\partial}{\partial x}\left( \frac{1}{|g|}\frac{\partial y}{\partial \eta} \right) = \frac{\partial \xi}{\partial x}\frac{\partial}{\partial \xi}\left( \frac{1}{|g|}\frac{\partial y}{\partial \eta} \right) + \frac{\partial \eta}{\partial x}\frac{\partial}{\partial \eta}\left( \frac{1}{|g|}\frac{\partial y}{\partial \eta} \right) \\
&= \frac{\partial \xi}{\partial x}\left( -\frac{1}{|g|^2}\frac{\partial g}{\partial \xi}\frac{\partial y}{\partial \eta} + \frac{1}{|g|}\frac{\partial^2 y}{\partial \xi \partial \eta} \right) + \frac{\partial \eta}{\partial x}\left( -\frac{1}{|g|^2}\frac{\partial g}{\partial \eta}\frac{\partial y}{\partial \eta} + \frac{1}{|g|}\frac{\partial^2 y}{\partial \eta^2} \right) \\
\frac{\partial^2 \eta}{\partial x^2} &= -\frac{\partial}{\partial x}\left( \frac{1}{|g|}\frac{\partial y}{\partial \xi} \right) = -\frac{\partial \xi}{\partial x}\frac{\partial}{\partial \xi}\left( \frac{1}{|g|}\frac{\partial y}{\partial \xi} \right) - \frac{\partial \eta}{\partial x}\frac{\partial}{\partial \eta}\left( \frac{1}{|g|}\frac{\partial y}{\partial \xi} \right) \\
&= -\frac{\partial \xi}{\partial x}\left( -\frac{1}{|g|^2}\frac{\partial g}{\partial \xi}\frac{\partial y}{\partial \xi} + \frac{1}{|g|}\frac{\partial^2 y}{\partial \xi^2} \right) - \frac{\partial \eta}{\partial x}\left( -\frac{1}{|g|^2}\frac{\partial g}{\partial \eta}\frac{\partial y}{\partial \xi} + \frac{1}{|g|}\frac{\partial^2 y}{\partial \xi \partial \eta} \right)
\end{aligned}, \tag{67a}$$

$$\begin{aligned}
\frac{\partial^2 \xi}{\partial x \partial y} &= \frac{\partial}{\partial y}\left( \frac{1}{|g|}\frac{\partial y}{\partial \eta} \right) = \frac{\partial \xi}{\partial y}\frac{\partial}{\partial \xi}\left( \frac{1}{|g|}\frac{\partial y}{\partial \eta} \right) + \frac{\partial \eta}{\partial y}\frac{\partial}{\partial \eta}\left( \frac{1}{|g|}\frac{\partial y}{\partial \eta} \right) \\
&= \frac{\partial \xi}{\partial y}\left( -\frac{1}{|g|^2}\frac{\partial g}{\partial \xi}\frac{\partial y}{\partial \eta} + \frac{1}{|g|}\frac{\partial^2 y}{\partial \xi \partial \eta} \right) + \frac{\partial \eta}{\partial y}\left( -\frac{1}{|g|^2}\frac{\partial g}{\partial \eta}\frac{\partial y}{\partial \eta} + \frac{1}{|g|}\frac{\partial^2 y}{\partial \eta^2} \right) \\
\frac{\partial^2 \eta}{\partial x \partial y} &= -\frac{\partial}{\partial y}\left( \frac{1}{|g|}\frac{\partial y}{\partial \xi} \right) = -\frac{\partial \xi}{\partial y}\frac{\partial}{\partial \xi}\left( \frac{1}{|g|}\frac{\partial y}{\partial \xi} \right) - \frac{\partial \eta}{\partial y}\frac{\partial}{\partial \eta}\left( \frac{1}{|g|}\frac{\partial y}{\partial \xi} \right) \\
&= -\frac{\partial \xi}{\partial y}\left( -\frac{1}{|g|^2}\frac{\partial g}{\partial \xi}\frac{\partial y}{\partial \xi} + \frac{1}{|g|}\frac{\partial^2 y}{\partial \xi^2} \right) - \frac{\partial \eta}{\partial y}\left( -\frac{1}{|g|^2}\frac{\partial g}{\partial \eta}\frac{\partial y}{\partial \xi} + \frac{1}{|g|}\frac{\partial^2 y}{\partial \xi \partial \eta} \right)
\end{aligned}, \tag{67b}$$

$$\begin{aligned}
\frac{\partial^2 \xi}{\partial y^2} &= -\frac{\partial}{\partial y}\left( \frac{1}{|g|}\frac{\partial x}{\partial \eta} \right) = -\frac{\partial \xi}{\partial y}\frac{\partial}{\partial \xi}\left( \frac{1}{|g|}\frac{\partial x}{\partial \eta} \right) - \frac{\partial \eta}{\partial y}\frac{\partial}{\partial \eta}\left( \frac{1}{|g|}\frac{\partial x}{\partial \eta} \right) \\
&= -\frac{\partial \xi}{\partial y}\left( -\frac{1}{|g|^2}\frac{\partial g}{\partial \xi}\frac{\partial x}{\partial \eta} + \frac{1}{|g|}\frac{\partial^2 x}{\partial \xi \partial \eta} \right) - \frac{\partial \eta}{\partial y}\left( -\frac{1}{|g|^2}\frac{\partial g}{\partial \eta}\frac{\partial x}{\partial \eta} + \frac{1}{|g|}\frac{\partial^2 x}{\partial \eta^2} \right) \\
\frac{\partial^2 \eta}{\partial y^2} &= \frac{\partial}{\partial y}\left( \frac{1}{|g|}\frac{\partial x}{\partial \xi} \right) = \frac{\partial \xi}{\partial y}\frac{\partial}{\partial \xi}\left( \frac{1}{|g|}\frac{\partial x}{\partial \xi} \right) + \frac{\partial \eta}{\partial y}\frac{\partial}{\partial \eta}\left( \frac{1}{|g|}\frac{\partial x}{\partial \xi} \right) \\
&= \frac{\partial \xi}{\partial y}\left( -\frac{1}{|g|^2}\frac{\partial g}{\partial \xi}\frac{\partial x}{\partial \xi} + \frac{1}{|g|}\frac{\partial^2 x}{\partial \xi^2} \right) + \frac{\partial \eta}{\partial y}\left( -\frac{1}{|g|^2}\frac{\partial g}{\partial \eta}\frac{\partial x}{\partial \xi} + \frac{1}{|g|}\frac{\partial^2 x}{\partial \xi \partial \eta} \right)
\end{aligned}. \tag{67c}$$

Differentiations of $\phi$ are given as follows:

$$\frac{\partial \phi}{\partial x} = \frac{\partial \xi}{\partial x}\frac{\partial \phi}{\partial \xi} + \frac{\partial \eta}{\partial x}\frac{\partial \phi}{\partial \eta}, \tag{68a}$$

$$\frac{\partial \phi}{\partial y} = \frac{\partial \xi}{\partial y}\frac{\partial \phi}{\partial \xi} + \frac{\partial \eta}{\partial y}\frac{\partial \phi}{\partial \eta}, \tag{68b}$$

$$\begin{aligned}
\frac{\partial^2 \phi}{\partial x^2} &= \frac{\partial}{\partial x}\left(\frac{\partial \xi}{\partial x}\frac{\partial \phi}{\partial \xi} + \frac{\partial \eta}{\partial x}\frac{\partial \phi}{\partial \eta}\right) = \frac{\partial^2 \xi}{\partial x^2}\frac{\partial \phi}{\partial \xi} + \frac{\partial \xi}{\partial x}\frac{\partial}{\partial x}\left(\frac{\partial \phi}{\partial \xi}\right) + \frac{\partial^2 \eta}{\partial x^2}\frac{\partial \phi}{\partial \eta} + \frac{\partial \eta}{\partial x}\frac{\partial}{\partial x}\left(\frac{\partial \phi}{\partial \eta}\right) \\
&= \frac{\partial^2 \xi}{\partial x^2}\frac{\partial \phi}{\partial \xi} + \frac{\partial \xi}{\partial x}\left(\frac{\partial \xi}{\partial x}\frac{\partial^2 \phi}{\partial \xi^2} + \frac{\partial \eta}{\partial x}\frac{\partial^2 \phi}{\partial \xi \partial \eta}\right) + \frac{\partial^2 \eta}{\partial x^2}\frac{\partial \phi}{\partial \eta} + \frac{\partial \eta}{\partial x}\left(\frac{\partial \xi}{\partial x}\frac{\partial^2 \phi}{\partial \xi \partial \eta} + \frac{\partial \eta}{\partial x}\frac{\partial^2 \phi}{\partial \eta^2}\right) \\
&= \frac{\partial^2 \xi}{\partial x^2}\frac{\partial \phi}{\partial \xi} + \left(\frac{\partial \xi}{\partial x}\right)^2\frac{\partial^2 \phi}{\partial \xi^2} + \frac{\partial \xi}{\partial x}\frac{\partial \eta}{\partial x}\frac{\partial^2 \phi}{\partial \xi \partial \eta} + \frac{\partial^2 \eta}{\partial x^2}\frac{\partial \phi}{\partial \eta} + \frac{\partial \xi}{\partial x}\frac{\partial \eta}{\partial x}\frac{\partial^2 \phi}{\partial \xi \partial \eta} + \left(\frac{\partial \eta}{\partial x}\right)^2\frac{\partial^2 \phi}{\partial \eta^2} \\
&= \left(\frac{\partial \xi}{\partial x}\right)^2\frac{\partial^2 \phi}{\partial \xi^2} + 2\frac{\partial \xi}{\partial x}\frac{\partial \eta}{\partial x}\frac{\partial^2 \phi}{\partial \xi \partial \eta} + \left(\frac{\partial \eta}{\partial x}\right)^2\frac{\partial^2 \phi}{\partial \eta^2} + \frac{\partial^2 \xi}{\partial x^2}\frac{\partial \phi}{\partial \xi} + \frac{\partial^2 \eta}{\partial x^2}\frac{\partial \phi}{\partial \eta},
\end{aligned} \tag{69a}$$

$$\begin{aligned}
\frac{\partial^2 \phi}{\partial x \partial y} &= \frac{\partial}{\partial x}\left(\frac{\partial \xi}{\partial y}\frac{\partial \phi}{\partial \xi} + \frac{\partial \eta}{\partial y}\frac{\partial \phi}{\partial \eta}\right) = \frac{\partial^2 \xi}{\partial x \partial y}\frac{\partial \phi}{\partial \xi} + \frac{\partial \xi}{\partial y}\frac{\partial}{\partial x}\left(\frac{\partial \phi}{\partial \xi}\right) + \frac{\partial^2 \eta}{\partial x \partial y}\frac{\partial \phi}{\partial \eta} + \frac{\partial \eta}{\partial y}\frac{\partial}{\partial x}\left(\frac{\partial \phi}{\partial \eta}\right) \\
&= \frac{\partial^2 \xi}{\partial x \partial y}\frac{\partial \phi}{\partial \xi} + \frac{\partial \xi}{\partial y}\left(\frac{\partial \xi}{\partial x}\frac{\partial^2 \phi}{\partial \xi^2} + \frac{\partial \eta}{\partial x}\frac{\partial^2 \phi}{\partial \xi \partial \eta}\right) + \frac{\partial^2 \eta}{\partial x \partial y}\frac{\partial \phi}{\partial \eta} + \frac{\partial \eta}{\partial y}\left(\frac{\partial \xi}{\partial x}\frac{\partial^2 \phi}{\partial \xi \partial \eta} + \frac{\partial \eta}{\partial x}\frac{\partial^2 \phi}{\partial \eta^2}\right) \\
&= \frac{\partial^2 \xi}{\partial x \partial y}\frac{\partial \phi}{\partial \xi} + \frac{\partial \xi}{\partial x}\frac{\partial \xi}{\partial y}\frac{\partial^2 \phi}{\partial \xi^2} + \frac{\partial \xi}{\partial y}\frac{\partial \eta}{\partial x}\frac{\partial^2 \phi}{\partial \xi \partial \eta} + \frac{\partial^2 \eta}{\partial x \partial y}\frac{\partial \phi}{\partial \eta} + \frac{\partial \xi}{\partial x}\frac{\partial \eta}{\partial y}\frac{\partial^2 \phi}{\partial \xi \partial \eta} + \frac{\partial \eta}{\partial x}\frac{\partial \eta}{\partial y}\frac{\partial^2 \phi}{\partial \eta^2} \\
&= \frac{\partial \xi}{\partial x}\frac{\partial \xi}{\partial y}\frac{\partial^2 \phi}{\partial \xi^2} + \left(\frac{\partial \xi}{\partial x}\frac{\partial \eta}{\partial y} + \frac{\partial \xi}{\partial y}\frac{\partial \eta}{\partial x}\right)\frac{\partial^2 \phi}{\partial \xi \partial \eta} + \frac{\partial \eta}{\partial x}\frac{\partial \eta}{\partial y}\frac{\partial^2 \phi}{\partial \eta^2} + \frac{\partial^2 \xi}{\partial x \partial y}\frac{\partial \phi}{\partial \xi} + \frac{\partial^2 \eta}{\partial x \partial y}\frac{\partial \phi}{\partial \eta},
\end{aligned} \tag{69b}$$

$$\begin{aligned}
\frac{\partial^2 \phi}{\partial y^2} &= \frac{\partial}{\partial y}\left(\frac{\partial \xi}{\partial y}\frac{\partial \phi}{\partial \xi} + \frac{\partial \eta}{\partial y}\frac{\partial \phi}{\partial \eta}\right) = \frac{\partial^2 \xi}{\partial y^2}\frac{\partial \phi}{\partial \xi} + \frac{\partial \xi}{\partial y}\frac{\partial}{\partial y}\left(\frac{\partial \phi}{\partial \xi}\right) + \frac{\partial^2 \eta}{\partial y^2}\frac{\partial \phi}{\partial \eta} + \frac{\partial \eta}{\partial y}\frac{\partial}{\partial y}\left(\frac{\partial \phi}{\partial \eta}\right) \\
&= \frac{\partial^2 \xi}{\partial y^2}\frac{\partial \phi}{\partial \xi} + \frac{\partial \xi}{\partial y}\left(\frac{\partial \xi}{\partial y}\frac{\partial^2 \phi}{\partial \xi^2} + \frac{\partial \eta}{\partial y}\frac{\partial^2 \phi}{\partial \xi \partial \eta}\right) + \frac{\partial^2 \eta}{\partial y^2}\frac{\partial \phi}{\partial \eta} + \frac{\partial \eta}{\partial y}\left(\frac{\partial \xi}{\partial y}\frac{\partial^2 \phi}{\partial \xi \partial \eta} + \frac{\partial \eta}{\partial y}\frac{\partial^2 \phi}{\partial \eta^2}\right) \\
&= \frac{\partial^2 \xi}{\partial y^2}\frac{\partial \phi}{\partial \xi} + \left(\frac{\partial \xi}{\partial y}\right)^2\frac{\partial^2 \phi}{\partial \xi^2} + \frac{\partial \xi}{\partial y}\frac{\partial \eta}{\partial y}\frac{\partial^2 \phi}{\partial \xi \partial \eta} + \frac{\partial^2 \eta}{\partial y^2}\frac{\partial \phi}{\partial \eta} + \frac{\partial \xi}{\partial y}\frac{\partial \eta}{\partial y}\frac{\partial^2 \phi}{\partial \xi \partial \eta} + \left(\frac{\partial \eta}{\partial y}\right)^2\frac{\partial^2 \phi}{\partial \eta^2} \\
&= \left(\frac{\partial \xi}{\partial y}\right)^2\frac{\partial^2 \phi}{\partial \xi^2} + 2\frac{\partial \xi}{\partial y}\frac{\partial \eta}{\partial y}\frac{\partial^2 \phi}{\partial \xi \partial \eta} + \left(\frac{\partial \eta}{\partial y}\right)^2\frac{\partial^2 \phi}{\partial \eta^2} + \frac{\partial^2 \xi}{\partial y^2}\frac{\partial \phi}{\partial \xi} + \frac{\partial^2 \eta}{\partial y^2}\frac{\partial \phi}{\partial \eta}.
\end{aligned} \tag{69c}$$

Substituting Eqs. (68) and (6 9) into Eq. (51), we obtain

$$\begin{aligned}
&\left[\left(\frac{\partial \xi}{\partial x}\right)^2 + \left(\frac{\partial \xi}{\partial y}\right)^2\right]\frac{\partial^2 \phi}{\partial \xi^2} + 2\left(\frac{\partial \xi}{\partial x}\frac{\partial \eta}{\partial x} + \frac{\partial \xi}{\partial y}\frac{\partial \eta}{\partial y}\right)\frac{\partial^2 \phi}{\partial \xi \partial \eta} + \left[\left(\frac{\partial \eta}{\partial x}\right)^2 + \left(\frac{\partial \eta}{\partial y}\right)^2\right]\frac{\partial^2 \phi}{\partial \eta^2} \\
&\quad + \left(\frac{\partial^2 \xi}{\partial x^2} + \frac{\partial^2 \xi}{\partial y^2}\right)\frac{\partial \phi}{\partial \xi} + \left(\frac{\partial^2 \eta}{\partial x^2} + \frac{\partial^2 \eta}{\partial y^2}\right)\frac{\partial \phi}{\partial \eta} = 0.
\end{aligned} \tag{70}$$

When $(x, y)$ is expressed explicitly as a function of $(\xi, \eta)$, the second derivatives of $(\xi, \eta)$ with respect to $(x, y)$ in Eq. (70) must be calculated using Eqs. (59) and (67).

### 5.1. Example 1: Introduction of polar coordinates

The boundary value problem is defined by Eqs. (50), (51) and (52) and is shown in Fig. 6(a). In the following, we consider the problem in the upper half plane $y > 0$ using the symmetry of the problem. For the purpose, we add condition on $y = 0$ :

$$\frac{\partial \phi}{\partial y} = 0 \quad \text{on} \quad a < |x| < +\infty, \quad y = 0. \tag{71}$$

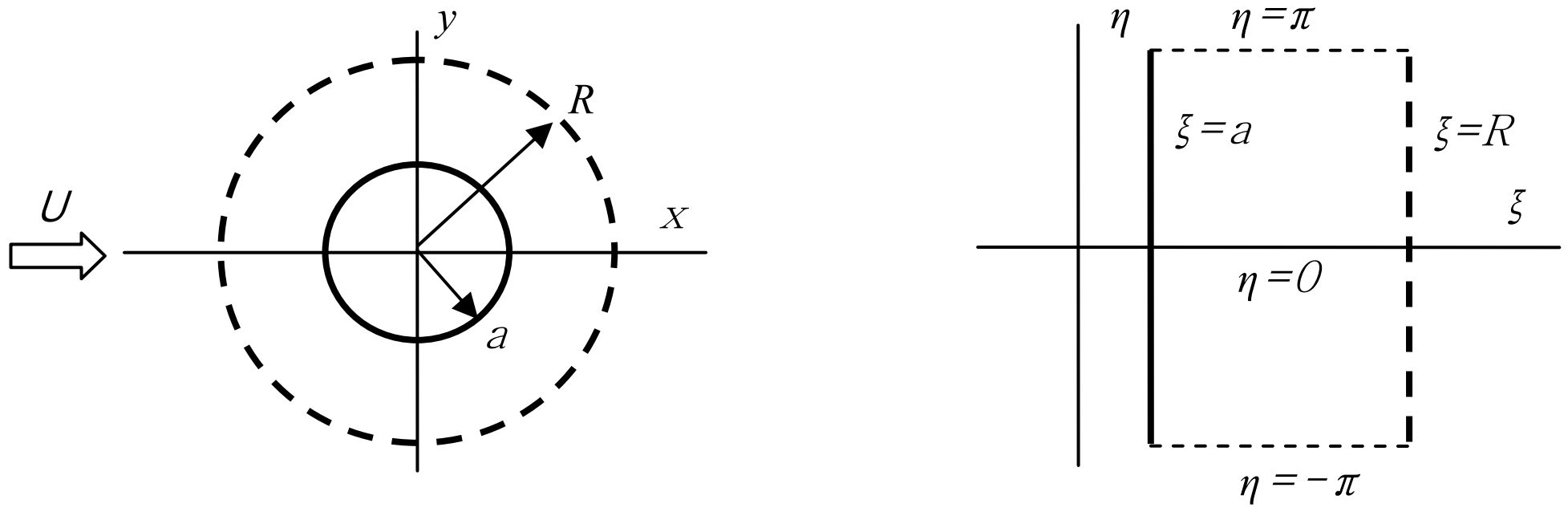

(a) $(x, y)$ plane (b) $(\xi, \eta)$ plane

Fig. 6. Polar coordinates

Firstly, we introduce the polar coordinates as shown in Fig. 6(b). The Cartesian coordinates $(x, y)$ is expressed explicitly using the polar coordinates $(\xi, \eta)$ as

$$x = \xi \cos\eta , \quad y = \xi \sin\eta , \tag{72}$$

$$\frac{\partial x}{\partial \xi} = \cos\eta , \quad \frac{\partial x}{\partial \eta} = -\xi \sin\eta , \quad \frac{\partial^2 x}{\partial \xi^2} = 0 , \quad \frac{\partial^2 x}{\partial \xi \partial \eta} = -\sin\eta , \quad \frac{\partial^2 x}{\partial \eta^2} = -\xi \cos\eta , \tag{73a}$$

$$\frac{\partial y}{\partial \xi} = \sin\eta , \quad \frac{\partial y}{\partial \eta} = \xi \cos\eta , \quad \frac{\partial^2 y}{\partial \xi^2} = 0 , \quad \frac{\partial^2 y}{\partial \xi \partial \eta} = \cos\eta , \quad \frac{\partial^2 y}{\partial \eta^2} = -\xi \sin\eta . \tag{73b}$$

In this case, the inverse transformation can be obtained easily, and the polar coordinates $(\xi, \eta)$ can be given by the Cartesian coordinates $(x, y)$ as

$$\xi = \sqrt{x^2 + y^2} , \quad \eta = \cos^{-1} \frac{x}{\sqrt{x^2 + y^2}} . \tag{74}$$

The derivatives of the polar coordinates with respect to the Cartesian coordinates could be obtained also be obtained easily. In this case, the coefficient of the differential equation (70) could be derived directly from Eq. (74). However, if we use Eqs. (72) and (73), the coefficients are obtained Eqs. (59), (60), (64) and (67). In the following numerical calculations, we adopt a method using Eqs. (72) and (73).

The exact solution is given below:

$$\Phi = U\xi \cos\eta + U \frac{a^2}{\xi} \cos\eta = Ux + Ua^2 \frac{x}{x^2 + y^2} , \tag{75}$$

$$u = \frac{\partial \Phi}{\partial x} = U + Ua^2 \frac{1}{x^2 + y^2} - Ua^2 \frac{2x^2}{(x^2 + y^2)^2} , \quad v = \frac{\partial \Phi}{\partial y} = -Ua^2 \frac{2xy}{(x^2 + y^2)^2} . \tag{76a,b}$$

In the numerical calculations, parameters are given as $R = 10$, $M = 40$, $N = 40$, $a = 2$ and $U = 1$. Numerical results are shown in Figs. 7, 8 and 9, where CC means Curvilinear Coordinates. The numerical results show a reasonable coincidence with the analytical solutions.

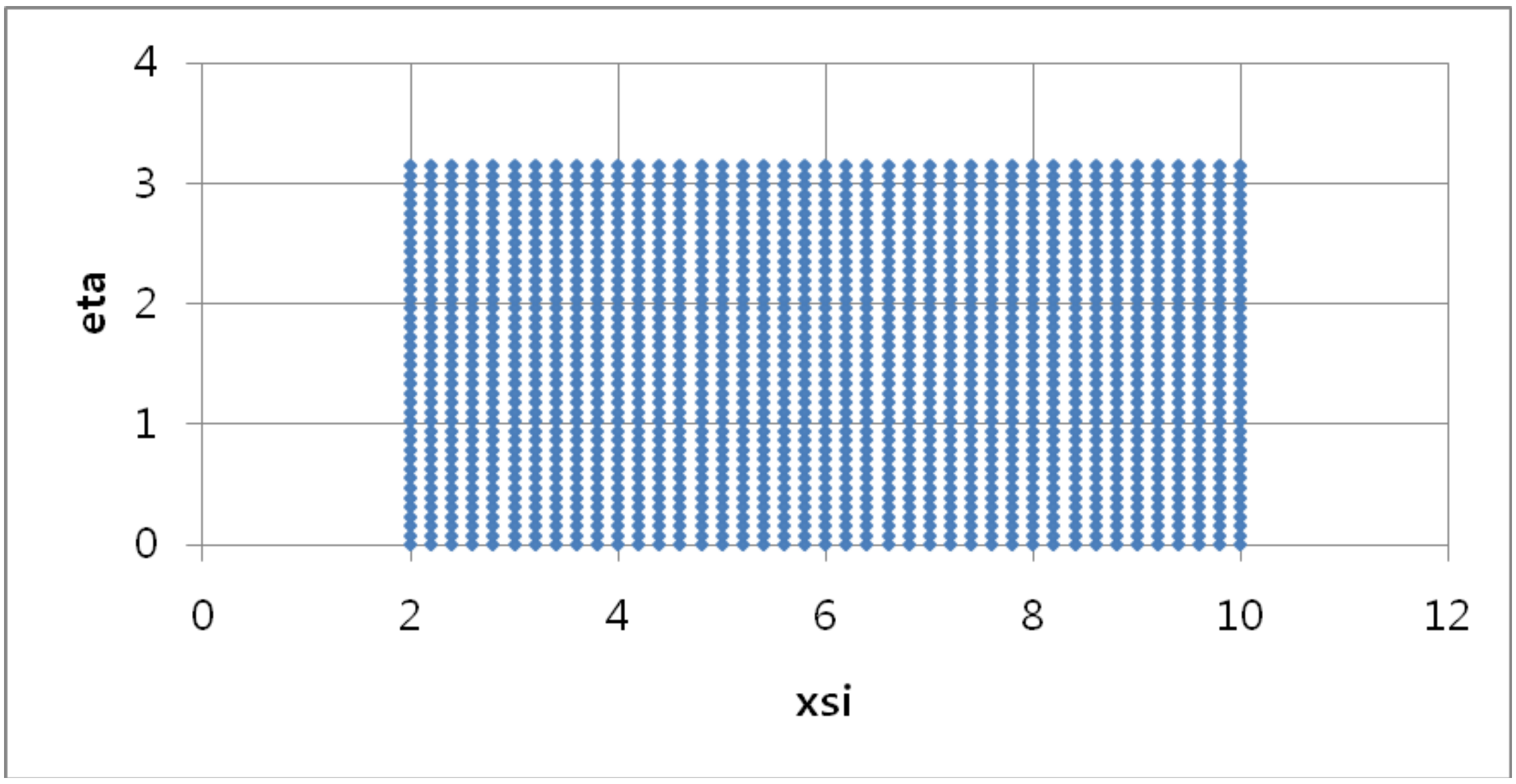

(a) Mapped plane

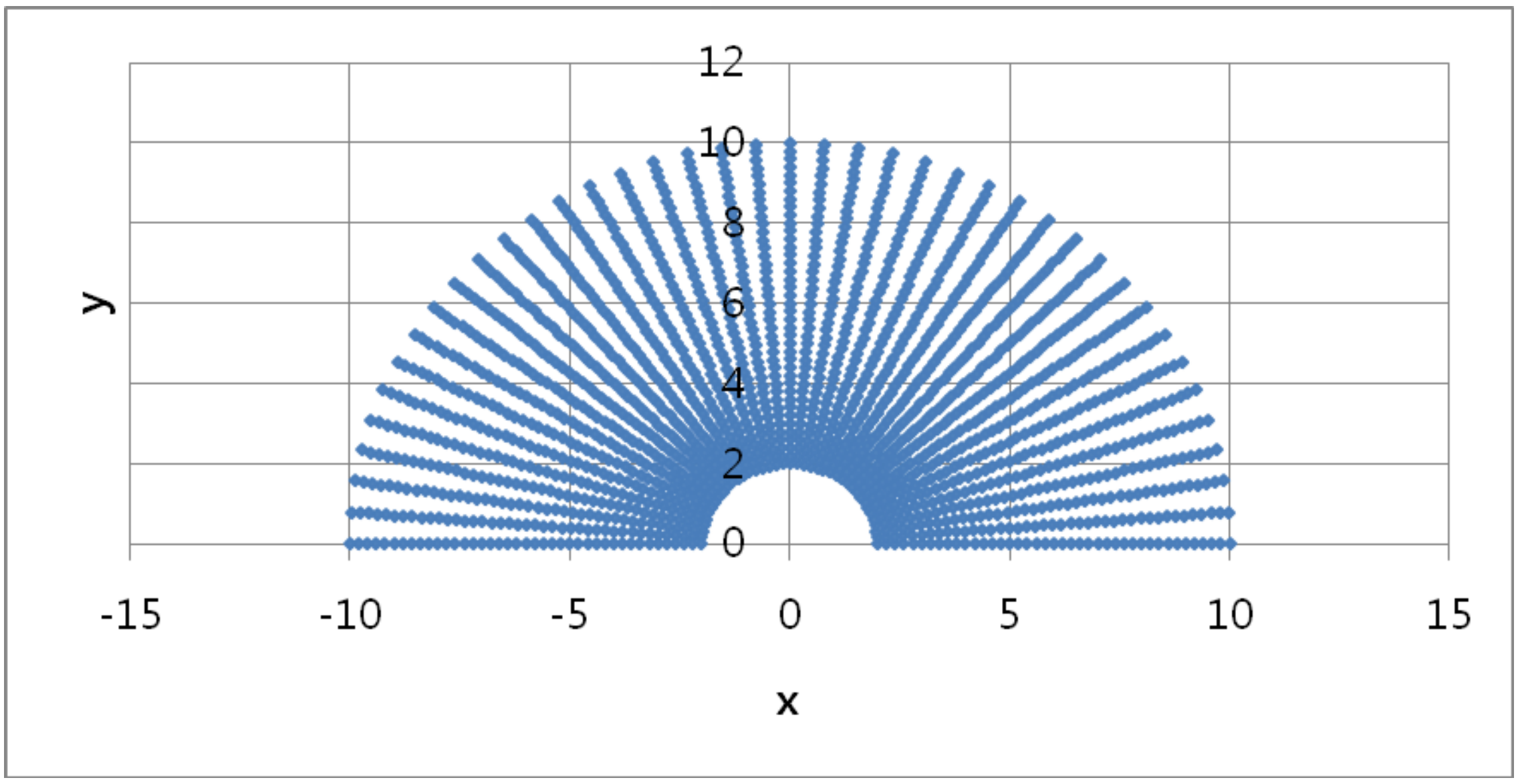

(b) Physical plane

Fig. 7. Coordinates transformation

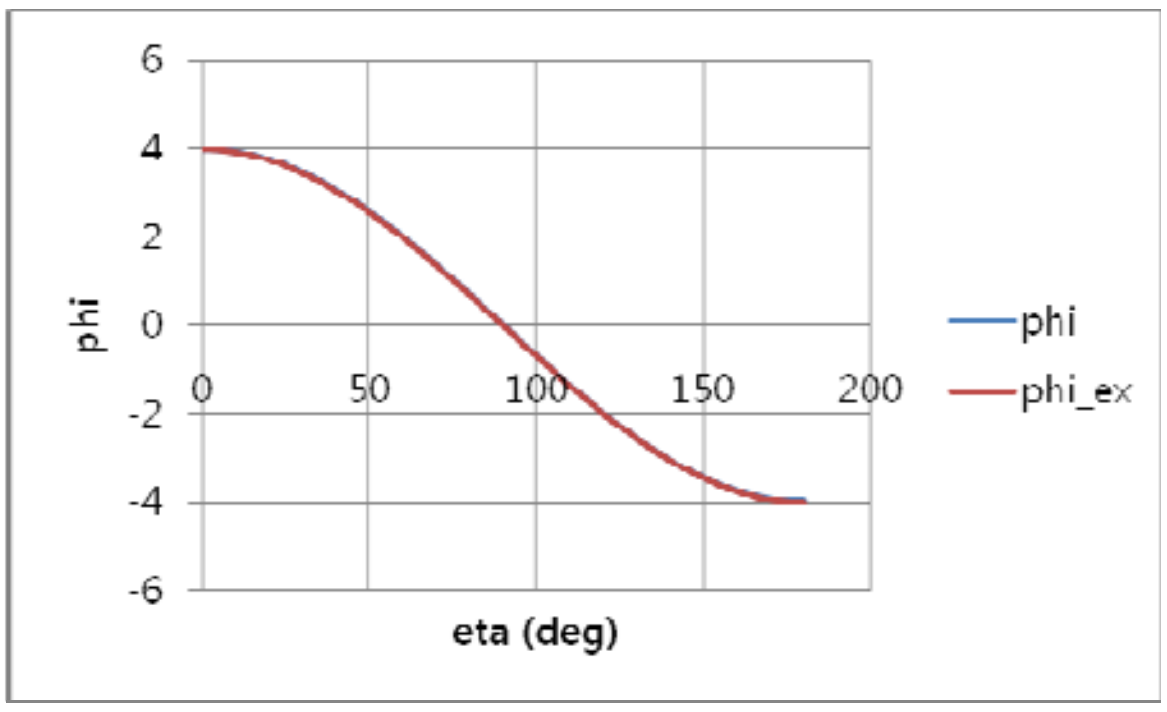

(a) Velocity potential on $r = a$

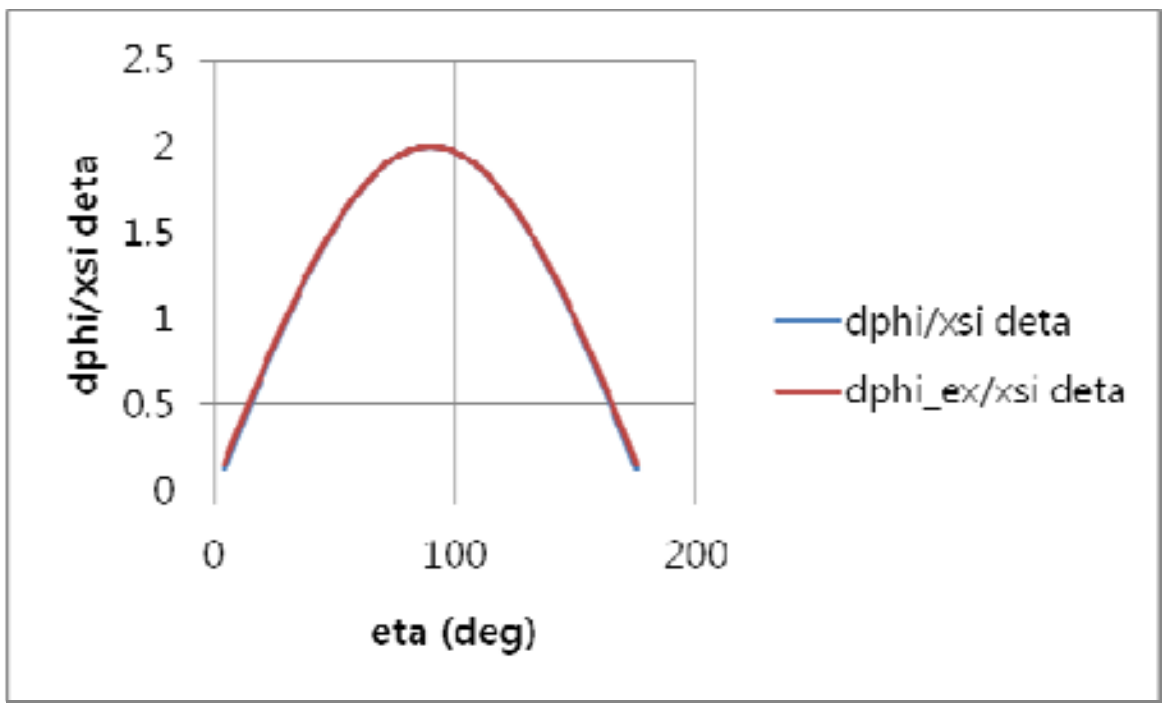

(b) Tangential velocity on $r = a$

Fig. 8. Comparison of velocity potential and tangential velocity with exact ones on $r = a$

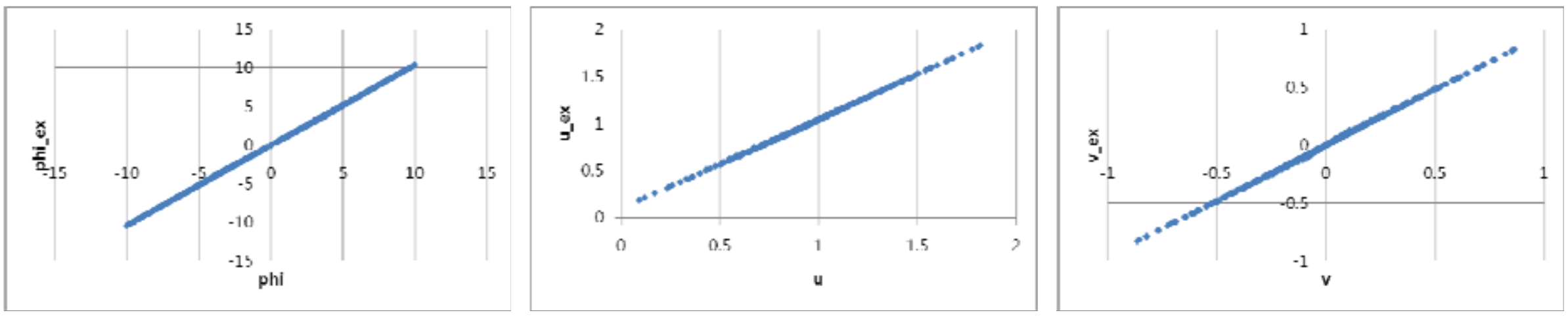

(a) Velocity potential phi (b) velocity u (c) velocity v

Fig. 9. Accuracy of the solution

**5.2. Example 2: Introduction of Joukowski transformation**

We solve the same problem as discussed in section 5.1 using a curvilinear coordinates obtained from Joukowski transformation:

$$\zeta = z + \frac{a^2}{z} . \tag{77}$$

where $a$ is the radius of the circular cylinder, and $z = x + iy$ and $\zeta = \xi + i\eta$ are the physical and mapped planes, respectively. The Joukowski transformation maps the circle $\sqrt{x^2 + y^2} = a$ in the $z$-plane as shown Fig. 10(a) into the line segment $|\xi| < 2a$, $\eta = 0$. The regular mesh in the $\zeta$-plane is shown in Fig. 10(b).

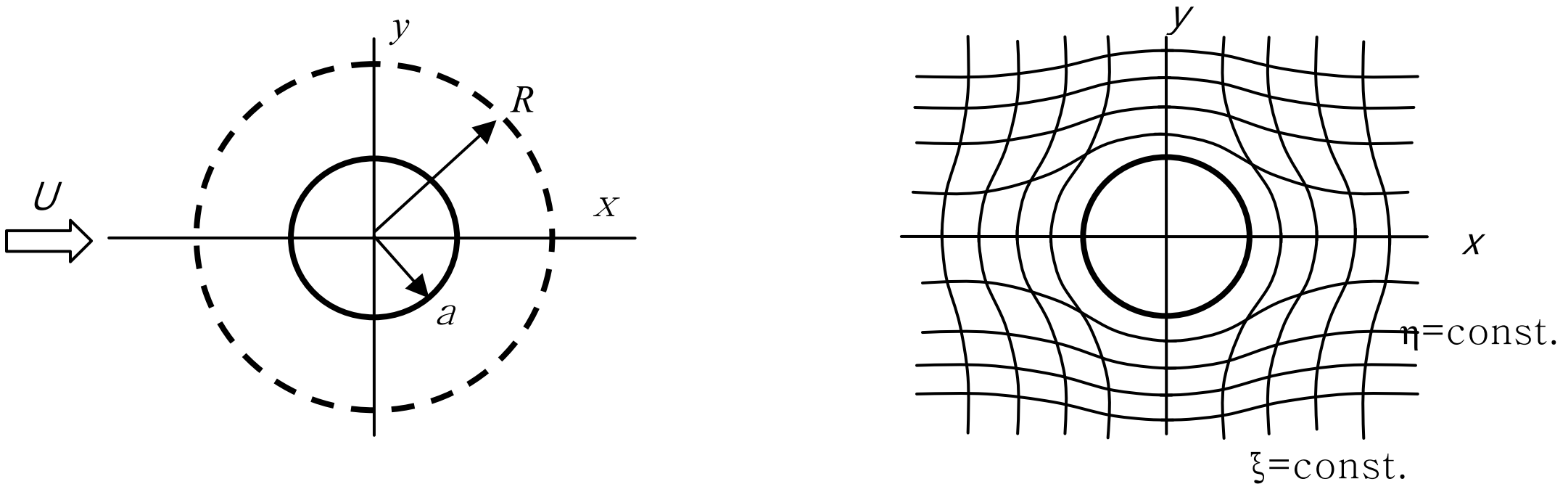

(a) $(x, y)$ plane (b) $\xi$ = const. and $\eta$ = const. lines)

Fig. 7. Jowkowski transformatiom

From Eq. (77), the transformation between $(\xi, \eta)$ and $(x, y)$ is given by

$$\xi = x + \frac{a^2 x}{x^2 + y^2}, \quad \eta = y - \frac{a^2 y}{x^2 + y^2} . \tag{78a,b}$$

Also from Eq. (77), the inverse transformation of Eq. (77) is given by

$$z = \frac{\zeta}{2} \pm \sqrt{\frac{\zeta^2}{4} - a^2} . \tag{79}$$

However, in the present paper, Newton-Raphson method is applied to Eq. (78) to obtain $(x, y)$ from $(\xi, \eta)$. The iterative procedure with the initial value $(x, y) = (\xi, \eta)$ is shown below:

$$\xi = x + \frac{a^2 x}{x^2 + y^2} + \left( 1 + \frac{a^2}{x^2 + y^2} - \frac{2a^2 x^2}{\left(x^2 + y^2\right)^2} \right) dx - \frac{2a^2 xy}{\left(x^2 + y^2\right)^2} dy , \tag{80a}$$

$$\eta = y - \frac{a^2 y}{x^2 + y^2} + \frac{2a^2 xy}{(x^2 + y^2)^2} dx + \left( 1 - \frac{a^2}{x^2 + y^2} + \frac{2a^2 y^2}{(x^2 + y^2)^2} \right) dy . \tag{80b}$$

From Eq. (77), the derivatives of $(\xi,\eta)$ with respect to $(x,y)$ are given by

$$\frac{\partial\xi}{\partial x}=1-\frac{a^2(x^2-y^2)}{(x^2+y^2)^2},\quad \frac{\partial\xi}{\partial y}=-\frac{2a^2xy}{(x^2+y^2)^2},$$
$$\frac{\partial^2\xi}{\partial x^2}=2\frac{a^2(x^3-3xy^2)}{(x^2+y^2)^3},\quad \frac{\partial^2\xi}{\partial x\partial y}=2\frac{a^2(3x^2y-y^3)}{(x^2+y^2)^3},\quad -\frac{\partial^2\xi}{\partial y^2}=2\frac{a^2(x^3-3xy^2)}{(x^2+y^2)^3}, \tag{81a}$$

$$\frac{\partial\eta}{\partial x}=\frac{2a^2xy}{(x^2+y^2)^2},\quad \frac{\partial\eta}{\partial y}=1-\frac{a^2(x^2-y^2)}{(x^2+y^2)^2},$$
$$\frac{\partial^2\eta}{\partial x^2}=-2\frac{a^2(3x^2y-y^3)}{(x^2+y^2)^3},\quad \frac{\partial^2\eta}{\partial x\partial y}=2\frac{a^2(x^3-3xy^2)}{(x^2+y^2)^3},\quad \frac{\partial^2\eta}{\partial y^2}=2\frac{a^2(3x^2y-y^3)}{(x^2+y^2)^3}. \tag{81b}$$

In the solution of the problem, the region was approximated by $|x|<L$, $0|<y<B$. The following boundary conditions were applied:

$$\frac{\partial\phi}{\partial y}=0 \quad \text{on} \quad -L<x<-a \ \ \& \ \ a<x\le L, \quad y=0; \tag{82a}$$

$$\frac{\partial\phi}{\partial x}n_x+\frac{\partial\phi}{\partial y}n_y=-Un_x \quad \text{on} \quad |x|<a, \quad y<\sqrt{x^2+y^2}; \tag{82b}$$

$$\frac{\partial\phi}{\partial x}=-U \quad \text{at} \quad |x|=a, \quad y=0; \qquad \phi=0 \quad \text{at} \quad x=-L, \quad y=0; \tag{82c}$$

$$\frac{\partial\phi}{\partial\xi}=0 \quad \text{on} \quad 0<\eta<B, \quad |\xi|=L; \tag{83}$$

$$\frac{\partial\phi}{\partial\eta}=0 \quad \text{on} \quad -L\le\xi\le L, \quad \eta=B. \tag{84}$$

In the numerical calculations, parameters are given as $L=10$, $B=10$, $M=40$, $N=40$, $a=2$ and $U=1$. Numerical results are shown in Figs. 10, 11 and 12. The numerical results show a reasonable agreement with the analytical solutions. However, the accuracy is lower than the results in section 5.1. The difference between Fig. 7(b) and Fig. 10(b) would explain the reason. The mesh in Fig. 10(b) is not fine enough at $x=\pm a$, $y=0$.

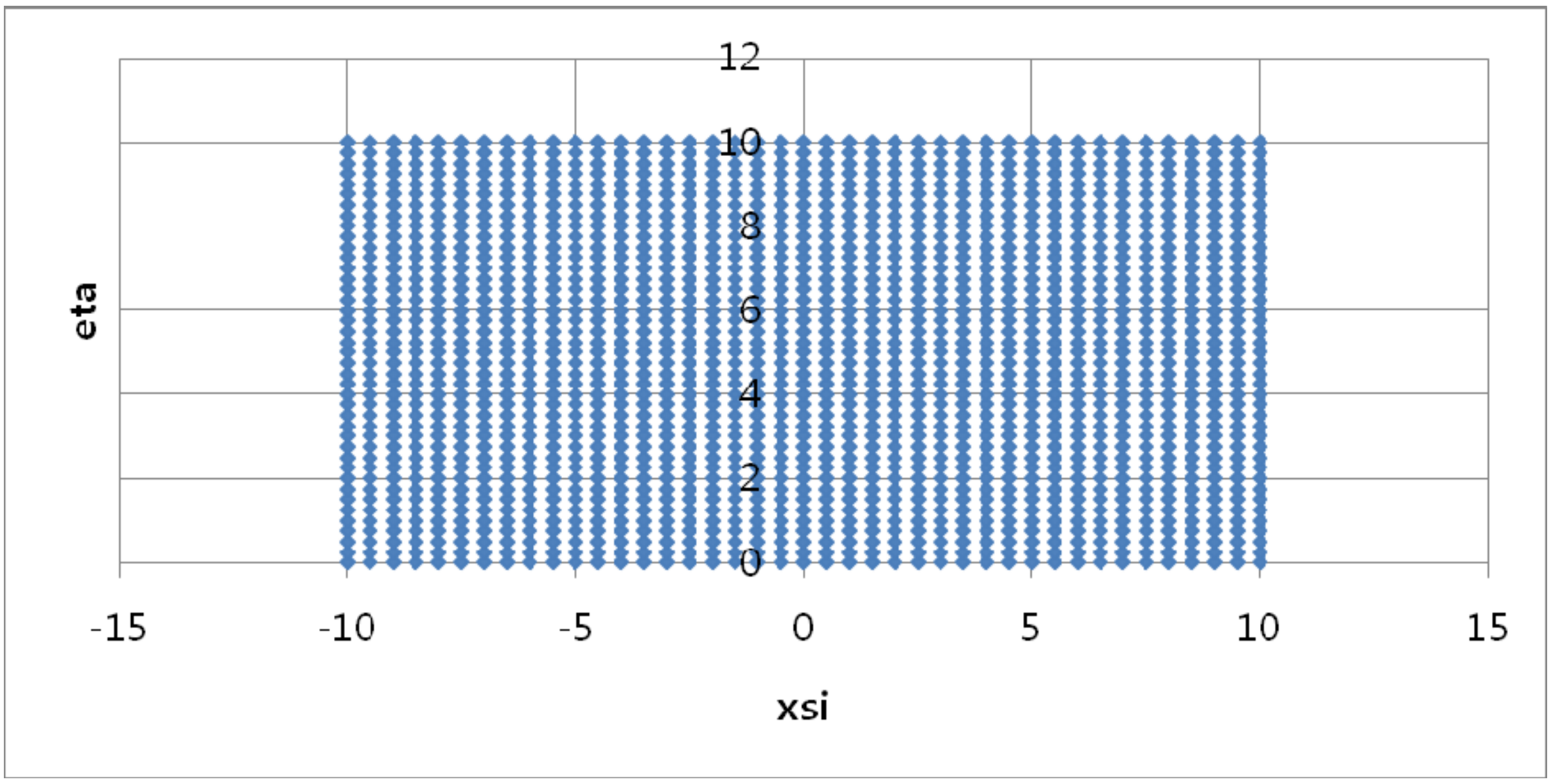

(a) Mapped plane

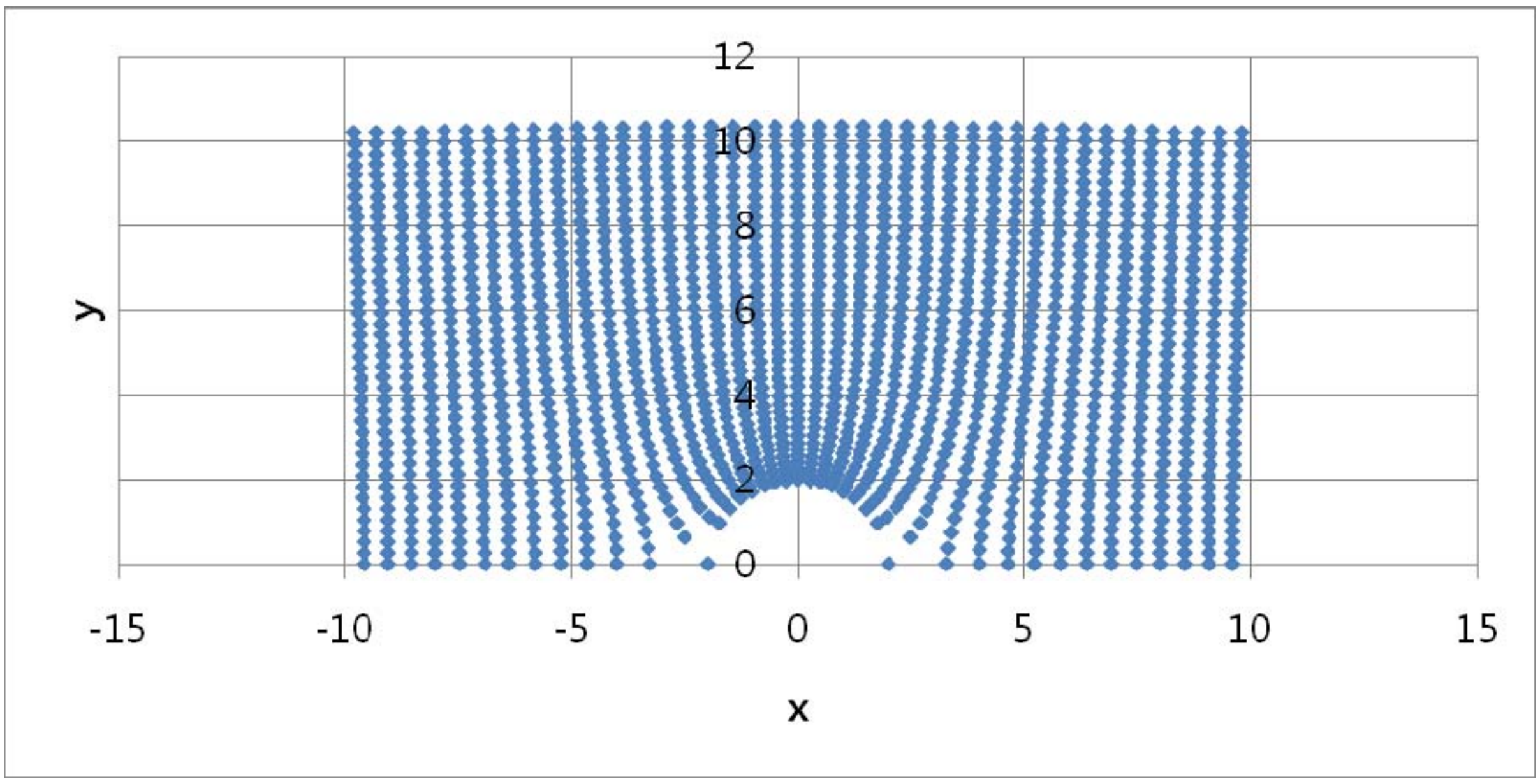

(b) Physical plane

Fig. 10. Coordinates transformation

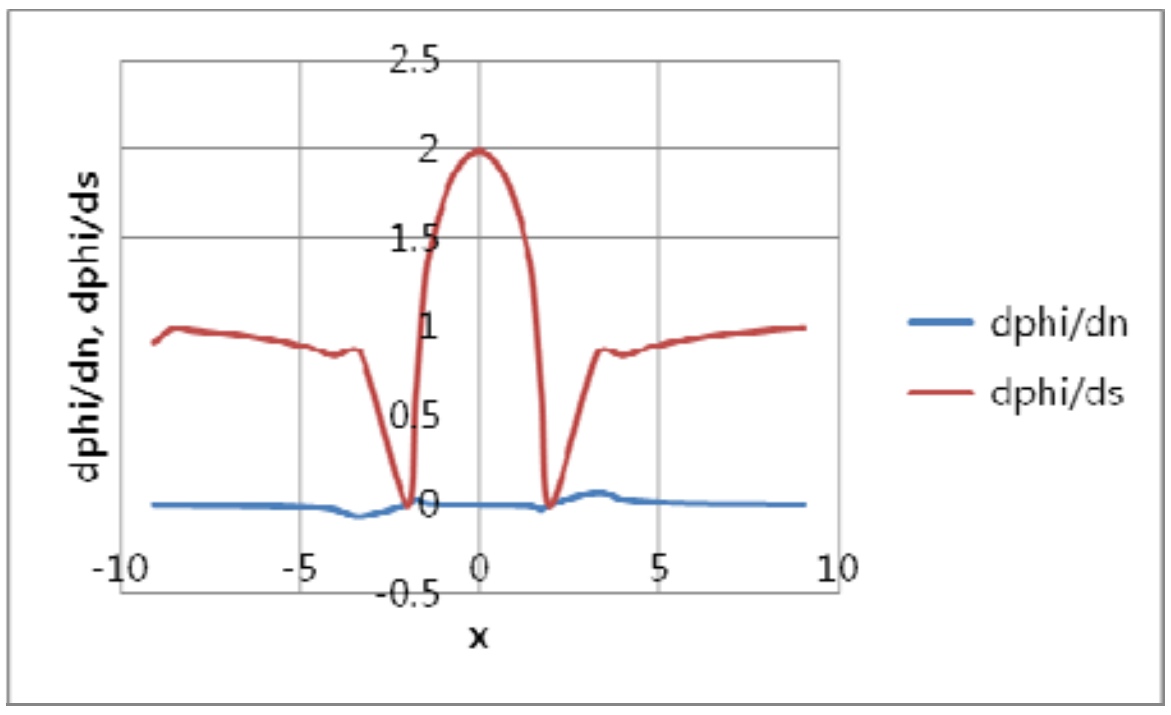

Fig.11. Normal velocity $u_n = d\phi/dn$ and tangential velocity $u_s = d\phi/ds$ on cylinder surface

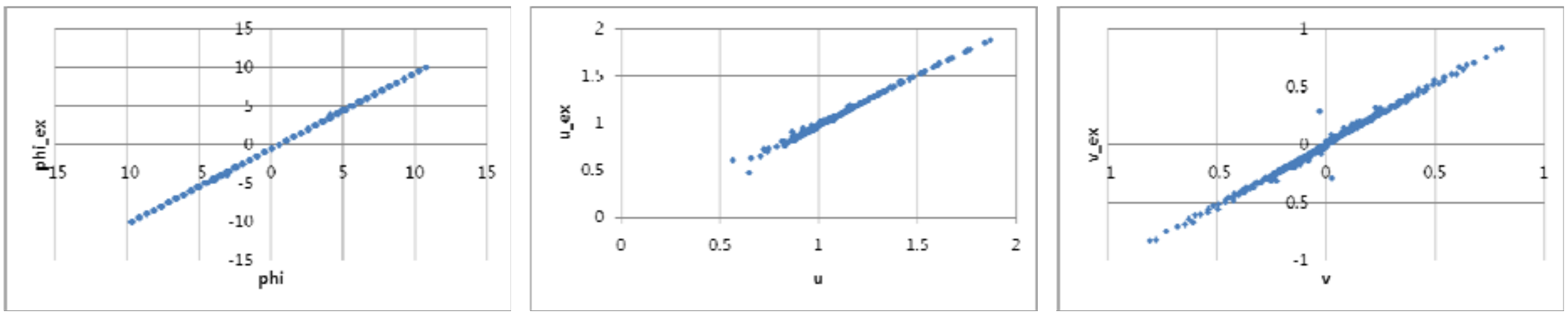

(a) Velocity potential phi (b) velocity u (c) velocity v

Fig. 12. Accuracy of the solution

**6. Conclusions**

Introduction of curvilinear coordinates into numerical computation was discussed, and the importance was shown. Although tensor analysis would be most suited for this kind of discussion, tensor notation can't be used in the present computer language. For example, the strict discrimination of upper and lower suffices is impossible in the present computer languages. In the present paper, a rather elementary approach more suited to write codes of computer programming is adapted. The theory is developed not specifically but generally. Hence, the theory could be applied widely.

In the present paper, we studied three problems, namely, fixed discontinuity, moving discontinuity and curved boundary. We can generate fine mesh in the neighborhood of fixed and moving discontinuities using Curvilinear Coordinates (CC). Furthermore, if we introduce curvilinear coordinates, we can transform a non-square region into square one and can use a regular mesh. Usually in numerical calculation, a curved boundary is approximated by a jagged or non-smooth boundary. The appropriateness and importance of introducing curvilinear coordinates into numerical analyses were shown reasonably through theoretical considerations and numerical verifications. An introduction of curvilinear coordinates could solve many important problems in the numerical analyses.

In the present paper, the transformation functions were given analytically. However, we would be able to define coordinates transformations using non-analytical methods. For example, we could draw a picture of mesh division and generate discrete data for the transformation. If we interpolate the discrete data using the proper base functions and interpolation using, for example, Least Square Method (LSM), we would be able to obtain a continuous and differentiable transformation function.